\documentclass{amsart}[12pt]

\usepackage{geometry} 

\usepackage{hyperref} 

\begin{document}

\title{ Evolution of starshaped hypersurfaces by general curvature functions }

\author{Ali Fardoun}
\address{Laboratoire de Math\'ematiques, UMR 6205 CNRS 
Universit\'e de Bretagne Occidentale 
6 Avenue Le Gorgeu, 29238 Brest Cedex 3   
France}
\email{Ali.Fardoun@univ-brest.fr}
\author{Rachid Regbaoui }
\address{Laboratoire de Math\'ematiques, UMR 6205 CNRS 
Universit\'e de Bretagne Occidentale 
6 Avenue Le Gorgeu, 29238 Brest Cedex 3   
France}
\email{Rachid.Regbaoui@univ-brest.fr}


\begin{abstract}
We consider the evolution of starshaped hypersurfaces in the Euclidean space by general curvature functions. Under appropriate conditions on the curvature function, we prove the global existence and convergence of the flow to a hypersurface of prescribed curvature.  \end{abstract} 
\maketitle
\section{Introduction}

 Let \ $M_0$ \ be a smooth closed  compact hypersurface in $\mathbb{R}^{n+1}$ ($n \geq 2$). We suppose that \ $M_0$ is starshaped with respect to a point, which we assume to be the origin of $\mathbb{R}^{n+1}$ for simplicity, and  in the rest of the paper all  starshaped hypersurfaces are   with respect to  the origin of  $\mathbb{R}^{n+1}$.  This means that for every point $P \in M_0$, we have $P \not\in  T_PM_0$, where $T_PM_0$ is the tangent space of $M_0$ at $P$.  If we let $\pi : M_0 \to \mathbb{S}^n$ to  be the projection on $\mathbb{S}^n$ defined by 
 
$$\pi(P) = {P \over |P|} \ , \ P \in M_0, $$

\noindent then one can prove that $M_0$ is starshaped if and only if $\pi$ is a diffeomorphism.  It follows that the  inverse diffeomorphism  $ X_0 := \pi^{-1} :    \mathbb{S}^n \to M_0 $   can be used as a parametrization of $M_0$.  The function  $\rho_0 : \mathbb{S}^n \to  \mathbb{R}^{+}$  defined by  \ $\rho_0(x) = |X_0(x)| $\ 
  is called the radial function of $M_0$.  Thus we have 
$$ X_0(x)  = \rho_0(x) x  ,   \ \ x \in \mathbb{S}^n  .   \eqno (1.1)$$

 \medskip
 
 From now on, we say that a smooth embedding  $ X :  \mathbb{S}^n \to \mathbb{R}^{n+1}$ is a starshaped embedding if $M := X(\mathbb{S}^n) $ is a starshaped hypersurface in $\mathbb{R}^{n+1}$, so by composing  by a smooth diffeomorphism of $\mathbb{S}^n$ if necessary, we may suppose that $X$ is of the form (1.1). 
 
 \medskip
 
 We consider the evolution  problem
$$\begin{cases} \partial_t X(t, x) = \Bigl(K\circ\kappa(X)(t,x) - f\circ X(t , x)\Bigr) \nu (t, x) \cr X(0, x) = X_0(x)  
\end{cases}
\eqno (1.2)$$
where  \ $X(t , . ): \mathbb{S}^n \to \mathbb{R}^{n+1}$ is a smooth starshaped embedding,  \ $ \nu $ is the outer unit normal vector field of the hypersurface $M_t := X(t, \mathbb{S}^n)$, \   $K$ is a suitable function of  the  principal curvatures vector  \ $ \kappa(X)  = (\kappa_1(X) , ...,  \kappa_n(X))$ of  $M_t$,  reffered as the curvature function, and \ $f :  \mathbb{R}^{n+1} \setminus \{ \ 0\ \} \to \mathbb{R}$ is a given smooth  function  referred as the prescribed function. We suppose that the function $K$ is expressed as an inverse function of the principal curvatures, that is 
$$K\circ\kappa(X) = {1 \over F\circ\kappa (X)}  = {1 \over F\circ(\kappa_1 (X), ..,  \kappa_n (X))}  \ , $$ 
where \ $F \in C^{\infty}(\Gamma) \cap C^0\left(\overline{\Gamma}\right)$ is a positive, symmetric function on an open, convex symmetric cone \ $\Gamma \subset \mathbb{R}^{n}$ with vertex at the origin, which contains the positive cone 
$$ \Gamma^{+} = \left\{ \ (\lambda_1,.., \lambda_n) \in  \mathbb{R}^{n} : \lambda_i > 0 \ \ \forall i \in [1,..,n] \  \right\}. $$

\noindent  This implies in particular that   
$$\Gamma \subset   \left\{ \ (\lambda_1,.., \lambda_n) \in  \mathbb{R}^{n} : \lambda_1 + \cdots + \lambda_n  > 0 \  \right\} . $$

\medskip

The function \  $F(\lambda)= F(\lambda_1,..,\lambda_n)$ is assumed to satisfy the following structure conditions 

$${\partial F \over \partial \lambda_i} > 0 \ \mathrm{on} \ \Gamma \ \ \  \forall i \in [1,..,n] \ \eqno (1.3)$$

$$ F \  \hbox{  is homogeneous of degree  $k > 0$ on} \  \Gamma  \   \hbox{ and }  \   F \equiv 0 \ \mathrm{on} \
 \partial\Gamma  \ \eqno (1.4)$$
 
$$ \log F  \    \   \hbox{is concave   on}  \  \Gamma . \eqno (1.5)$$ 
 \medskip
 
\noindent By scaling, we may suppose 
$$F(1,..,1) = 1 . \eqno  (1.6) $$  

\bigskip

 The above conditions on $F$ are usually assumed in the study of fully nonlinear partial differential equations.  Condition (1.3) ensures that the system (1.2) is parabolic, which  is  of great  importance in proving short time existence of solutions.  The other conditions will be used to control   the $C^1$ and $C^2$-norms of  solutions. Some examples of suitable curvature functions satisfying (1.3)-(1.6)   are   
   $$F(\lambda_1,..,\lambda_n)= \binom{k}{n}^{-1} S_k(\lambda_1,...,\lambda_n) = \binom{k}{n}^{-1} \sum_{1 \le i_1 < \cdots < i_k \le n }\lambda_{i_1} \cdots \lambda_{i_k} $$
   the $k$-th elementary symmetric functions normalised so that \ $F(1,.., 1)= 1$. In this case we take \ $\Gamma$ \ to be the component of the set where \ $S_k$ is positive which contains the positive cone. Thus we obtain the mean curvature when \ $k=1$ and the Gauss curvature when $k=n$. Other examples of curvature functions are
  $$F(\lambda_1,..,\lambda_n)= \binom{k}{n} \left(S_k ( \lambda_1^{-1},..., \lambda_n^{-1} )\right )^{-1}. $$
  In this case, we take \ $\Gamma =  \Gamma^{+}$. A particular case of interest in the previous example  is the harmonic mean curvature when \ $k = 1$. 
  
  \medskip
  
   Finally, we notice that if a function $F$ satisfies conditions (1.3)-(1.6) above, then for any $\alpha > 0$, the function $F^{\alpha}$ satifies the same conditions where $k$ is replaced by $\alpha k$.  This invariance property  is due to the fact that the convexity condition (1.5) concerns $\log F$ but not $F$.

  \medskip

  When  the prescribed function \ $f \equiv 0$, problem (1.2) has been studied by J. Urbas \cite{jU90} and independently by  C.  Gerhardt \cite{cG90}, assuming that the curvature  function \ $F$ satisfies (1.3)-(1.6)  with $k = 1$ and that  $F$ is concave instead of $\log F$ concave.  They showed the existence of a global solution on \ $[0, + \infty)$, and for the convergence at infinity, they prove  that if \ $\tilde M_t$ \ is the hypersurface parametrized by \ $\tilde X(t , . )= e^{-t}X(t , . )$,  then \ $\tilde M_t$ \ converges to a sphere in the \ $C^{\infty}$ topology as \ $t \to + \infty$.  See also the recent work of  C.  Gerhardt \cite{cG12} where he considers curvature functions $F$ with  homogeneity degree $ 0 < k \not=1$ and $f \equiv 0$. There is an extensive litterature on curvature evolution equation like (1.2) when $f\equiv 0$. We refer the reader to \cite{bA94}, \cite{bA07}, \cite{gH96},  \cite{kS98}, \cite{kS02}, \cite{jU91} and the references therein. 
  
  \medskip
  
   In this paper, we study the global existence and convergence for equation (1.2)  assuming that  $F$ satisfies the structure conditions (1.3)-(1.6), and the prescribed function $f :  \mathbb{R}^{n+1} \setminus \{ \ 0\ \} \to \mathbb{R}^{+}$ is a smooth function satisfying 
     $${\partial \over \partial \rho}\left( \rho^{-k} f(X) \right) > 0 , \   X \in   \mathbb{R}^{n+1} \setminus \{ \ 0\ \}  \eqno (1.7)$$
      where \ $ \rho = |X|$. We will  also assume that  there exist two positive  real numbers $r_1 \le r_2$  such that 
    $$\begin{cases}f(X) \le r_1^k \   \   \hbox{if} \  |X| = r_1 \cr $$
    f(X) \ge r_2^k  \  \   \hbox{if}  \  |X|= r_2  . \end{cases} \eqno (1.8)$$
      
   \noindent  These assumptions were made by   L. Caffarelli, L. Nirenberg and J. Spruck \cite{lC} for the existence by elliptic methods of starshaped embedding $X$  whose   ${1\over F}$-curvature is equal to $f$, i.e, statisyfing the equation :
    $$ {1 \over F(\kappa (X))}  = f(X)  . \eqno (1.9) $$
    
    Our main result in this paper is that  conditions (1.7)-(1.8)  on the prescribed function $f$ are also sufficient to study the evolution problem (1.2). Moreover the solution of such flow converges to a smooth starshaped embedding  satisfying (1.9).  Our first  result concerns the case where the homogeneity degree $k$ of $F$ satisfies $ 0 < k \le 1$.  We have 
  
  \medskip
  
  \newtheorem{theo}{Theorem}[section] \begin{theo} Let \ $F \in C^{\infty}(\Gamma) \cap C^0\left(\overline{\Gamma}\right)$ be a positive symmetric function satisfying  conditions $(1.3)$-$(1.6)$ such that the  homogeneity degree $k$ of $F$ satisfies $ 0 < k \le 1$, and  let  $f  \in C^{\infty}\left(\mathbb{R}^{n+1} \setminus \{  0  \}\right)$ be a positive smooth function    satisfying $(1.7)$-$(1.8)$.  Let $M_0$ a closed compact starshaped hypersurface in $\mathbb{R}^{n+1},$  paramatrized by an embedding   $X_0  : \mathbb{S}^n \to
 \mathbb{R}^{n+1}$ of the form $ (1.1)$ such that 
 $$  \kappa(X_0) \in \Gamma   \   \    \hbox{and} \    \   { 1 \over F(\kappa(X_0))} - f(X_0) \ge 0 \   .   \eqno (1.10) $$

     \medskip
     
  \noindent Then  the evolution problem $(1.2)$ admits a unique smooth  solution $X(t , . )$ defined on $[0 , + \infty )$ such that,  for every  \  $ t \in [0, +\infty),  \  X(t , . )  :    \mathbb{S}^n \to
 \mathbb{R}^{n+1} $  is   a starshaped embedding  satisfying   $  \kappa(X(t , . ))  \in \Gamma$.    Moreover,  $X(t, .)$ converges in  $C^{\infty}(\mathbb{S}^n, \mathbb{R}^{n+1})$ to a starshaped embedding $ X_{\infty} : \mathbb{S}^n \to \mathbb{R}^{n+1}$ as $t \to + \infty$, satisfying
 $${1 \over F(\kappa(X_{\infty}))} = f(X_{\infty}) , $$
 and for any  $m \in \mathbb{N}$, $t \in [0, + \infty )$, we have 
 $$\| X(t , .) - X_{\infty}\|_{C^m(\mathbb{S}^n, \ \mathbb{R}^{n+1})} \le C_me^{-\lambda_m t} ,  \eqno (1.11) $$
 where $C_m$ and $\lambda_m$ are positive constants depending only on $m, f, F, r_1, r_2$ and $X_0$. 
 
   \end{theo}

  \bigskip

\newtheorem{rem}{Remark}[section]
 \begin{rem} There are many  starshaped embeddings  $X_0  : \mathbb{S}^n\to \mathbb{R}^{n+1}$ satisfying condition (1.10)  in Theorem 1.1. Indeed, it suffices to  take  $X_0(x) = rx , x \in \mathbb{S}^n$,  where $r$ is any positive constant such that $0 < r  \le r_1$,  with $r_1$    as in (1.8). Using conditions (1.7)-(1.8), it is easy to see that (1.10) is satsified.  \end{rem}

\medskip
 
 As a  consequence of Theorem 1.1, we recover the existence result for Weingarten hypersurfaces of Cafarelli-Nirenberg-Spruck \cite{lC} stated above.  Moreover,  we prove  the uniqueness of starshaped solutions of (1.9). Namely   we have : 

\medskip

\medskip

\newtheorem{cor}{Corollary}[section] \begin{cor}
  Let \ $F \in C^{\infty}(\Gamma) \cap C^0\left(\overline{\Gamma}\right)$ be a positive symmetric function satisfying $(1.3)$-$(1.6)$, and let  $f $ be a  smooth positive  function  satisfying  $(1.7)$-$(1.8)$.   Then there exists a smooth starshaped embedding $X :  \mathbb{S}^n \to
 \mathbb{R}^{n+1}$ such that $  \kappa(X) \in \Gamma $, and satisfying  
  $$ {1 \over F(\kappa(X))} =  f(X)  \  . \eqno (1.12)$$
  
  \noindent Moreover, $X$ is the unique starshaped solution of $(1.12)$  with $  \kappa(X) \in \Gamma $. 
\end{cor}

\medskip

 When the homogeneity degree $k$ of the curvature function $F$ satisfies $ k >1$, we need additional conditions on the initial embedding $X_0$. More precisely, we have

\medskip

\begin{theo} Let \ $F \in C^{\infty}(\Gamma) \cap C^0\left(\overline{\Gamma}\right)$ be a positive symmetric function satisfying  conditions $(1.3)$-$(1.6)$ such that the  homogeneity degree $k$ of $F$ satisfies $  k > 1$, and  let  $f  \in C^{\infty}\left(\mathbb{R}^{n+1} \setminus \{  0  \}\right)$ be a positive smooth function    satisfying $(1.7)$-$(1.8)$.  Let $M_0$ be a closed compact starshaped hypersurface in $\mathbb{R}^{n+1},$  paramatrized by an embedding   $X_0  : \mathbb{S}^n \to
 \mathbb{R}^{n+1}$ of the form $ (1.1)$ such that 
  $$  \kappa(X_0) \in \Gamma   \  \   \hbox{and} \  \   
 0 \le -\left({1 \over F(\kappa(X_0))} - f(X_0)\right){|\nabla X_0| \over |X_0|}  \le  {k R_1 \over (k+1) R_2} \min_{R_1 \le |Y| \le R_2}f(Y) \ ,  \eqno (1.13) $$
 
\noindent  where  
 $$ R_1 = \min\left(   r_1 , \min_{x \in \mathbb{S}^n}|X_0(x)|\right) \ , \  R_2 =  \max\left(   r_2 , \max_{x \in \mathbb{S}^n}|X_0(x)|\right)$$
 
\noindent  and $r_1, r_2$ are as in (1.8). Then  the evolution problem $(1.2)$ admits a unique smooth  solution $X(t , . )$ defined on $[0 , + \infty )$ such that,  for every  \  $ t \in [0, +\infty),  \  X(t , . )  :    \mathbb{S}^n \to
 \mathbb{R}^{n+1} $  is   a starshaped embedding  satisfying   $  \kappa(X(t , . ))  \in \Gamma$.    Moreover,  $X(t, .)$ converges in  $C^{\infty}(\mathbb{S}^n, \mathbb{R}^{n+1})$ to a starshaped embedding $ X_{\infty} : \mathbb{S}^n \to \mathbb{R}^{n+1}$ as $t \to + \infty$, satisfying
 $${1 \over F(\kappa(X_{\infty}))} = f(X_{\infty}) , $$
 and for any  $m \in \mathbb{N}$, $t \in [0, + \infty )$, we have 
 $$\| X(t , .) - X_{\infty}\|_{C^m(\mathbb{S}^n, \ \mathbb{R}^{n+1})} \le C_me^{-\lambda_m t} ,  \eqno (1.14) $$
 where $C_m$ and $\lambda_m$ are positive constants depending only on $m, f, F, r_1, r_2$ and $X_0$. 
 
   \end{theo}

\medskip

\medskip

 \begin{rem} There are many smooth  starshaped embeddings  $X_0  : \mathbb{S}^n\to \mathbb{R}^{n+1}$ satisfying condition (1.13)  in Theorem 1.2. Indeed,   by applying Corollary 1.1 to the functions $F^{1/k}, f^{1/k}$ instead of $F, f $  $($ as it can easily be seen,  conditions (1.3)-(1.6)  and (1.7)-(1.8) are  still satisfied  with a new homgeneñity degree $k = 1$ for $F^{1/k} )$, then we get a smooth starshaped  embedding $X: \mathbb{S}^n \to \mathbb{R}^{n+1}$ satisfying : 
  $${1 \over F(\kappa(X))} = f(X).  $$
If we  take $X_0 = r X$ , where $r$ is any positive  constant such that $r  \in [1, 1 + \varepsilon)$, with  $\varepsilon> 0$  small enough, then it is not difficult to see, by using condition (1.7)-(1.8), that $X_0$ satisfies condition (1.13) in Theorem 1.2. 
      \end{rem}


\section{Preliminaries}

In this section, we recall some expressions for the relevant geometric quantities of smooth compact starshaped hypersurfaces \ $M \subset \mathbb{R}^{n+1}$. As we saw in the previous section,   there is a smooth  embedding   $X : \mathbb{S}^n \to \mathbb{R}^{n+1}$ parametrizing  $M$,  which is of the form  
 $$X(x) = \rho(x)x, \ x \in \mathbb{S}^n . $$
 For any local orthonormal frame $\{ e_1, ..., e_n \} $ on  $\mathbb{S}^n$ (endowed with its standard metric), covariant differentiation with respect to $e_i$  will be denoted by $\nabla_i$,  $\nabla_{ij}, \nabla_{ijk}, ...$, and we let  $ \nabla$ be the gradient on  $\mathbb{S}^n$.  Then in terms of the radial function  $\rho$, the metric \ $ g = \left[g_{ij}\right]$  induced by $X$ and its inverse \ $g^{-1}= \left[g^{ij}\right] $   are given by 
$$g_{ij}= \langle \nabla_iX, \nabla_jX\rangle = \rho^2 \delta_{ij}+ \nabla_i\rho \nabla_j\rho, \ \ \ \ \ g^{ij}= \rho^{-2}\left(\delta_{ij}-{\nabla_i\rho \nabla_j\rho \over \rho^2 + | \nabla \rho |^2 }\right) , \eqno (2.1)$$
where \ $\langle \ , \  \rangle$ is the standard metric on \ $\mathbb{R}^{n+1}$, and $\delta_{ij}$ are Kronecker symbols.   The unit outer normal to $M$ is 
$$
\nu = {\rho x - \nabla \rho \over \sqrt{\rho^2 + | \nabla \rho |^2 }}  \eqno (2.2)$$ 
and the the second fundamental form of \ $M$ is given by 
$$
h_{ij} = - \langle \nabla_{ij} X, \nu \rangle = \left ( \rho^2 + | \nabla \rho |^2  \right)^{-{1 \over 2}} (\rho^2\delta_{ij} + 2 \nabla_i\rho\nabla_j \rho - \rho \nabla_{ij}\rho), \eqno(2.3)$$

\medskip

\noindent The principal curvatures of \ $M$ are the eigenvalues of the second fundamental form with respect to the induced metric $g$. Thus, $\lambda $ is a principal curvature if 
$$
\mathrm{det}[h_{ij} - \lambda g_{ij}] = 0 ,$$
or equivalently
$$
\mathrm{det}[a_{ij} - \lambda \delta_{ij}] = 0 ,
$$
where the symmetric matrix \ $ [a_{ij}]$  is given by 
$$
 [a_{ij}] = [g^{ij} ]^{{1 \over 2}} [h_{ij} ] [g^{ij} ]^{{1 \over 2}} \eqno (2.4)$$
and where \ $ [ g^{ij} ]^{{1 \over 2}}$ is the positive square root of \ $ [ g^{ij} ]$ which is given by
$$
 [g^{ij} ]^{{1 \over 2}}= \rho^{-1} \left[ \delta_{ij} - {\nabla_i\rho\nabla_j \rho \over \sqrt{\rho^2 + | \nabla \rho |^2 }(\rho  + \sqrt {\rho^2 + | \nabla \rho |^2 })} \right].
 \eqno (2.5) $$
 
 \medskip
 
 Let us now make some remarks about the curvature function $F$.  Since $F$ is symmetric,  it is well known that $F$ can be seen as a smooth function on  the set of real symmetric $n\times n$ matrices $[a_{ij}]$.   More precisely, we have 
 $$F \in C^{\infty}(M(\Gamma)) \cap C^0(\overline{M(\Gamma}))$$
 where $ M(\Gamma)$ is the convex cone of real symmetric $n\times n$ matrices with eigenvalue vector in the cone $\Gamma$. One can also prove that conditions (1.3)-(1.6) are also valid when $F$ is seen as function on $M(\Gamma)$. We have 

$$ \left[ F_{ij}\right] \ \hbox {is  positive definite  on } \ M(\Gamma) \ ,  \eqno (2.6)$$
 where $F_{ij} = {\partial F \over \partial a_{ij}}$.
 
$$ F \  \hbox{  is homogeneous of degree  $k > 0$ on} \  M(\Gamma)   \   \hbox{ and }  \   F \equiv 0 \ \mathrm{on} \
 \partial M(\Gamma)  \ \eqno (2.7)$$
 
$$ \log F  \    \   \hbox{is concave   on}  \  M(\Gamma) . \eqno (2.8)$$

 $$F(\delta_{ij}) = 1 . \eqno  (2.9) $$ 

\medskip

 \noindent We note here that a smooth function  $G $ on $M(\Gamma)$  is concave if
 $$\sum_{i,j=1}^n\sum_{k,l=1}^nG_{ij, kl}\  \eta_{ij}\eta_{kl} \le 0 \hspace{3mm} \text{on} \ M(\Gamma) $$
 for all real symetric $n\times n$ matrices $(\eta_{ij})$, where 
 $$ G_{ij, kl} = {\partial^2 G \over \partial a_{kl} \partial a_{ij} }.$$

  Now, we will show that equation (1.2) is equivalent to an evolution equation depending on the radial function $\rho$. We proceed as in \cite{jU90}, first suppose that \ $X(t , . )$ is a solution of (1.2) such that for each \ $t \in [0, + \infty)$, $X(t , . )$ is an embedding of a smooth closed compact hypersurface \ $M_t$ in \ $\mathbb{R}^{n+1}$, which is starshaped with respect to the origin and  such that the vector of its principal curvatures $\kappa = ( \kappa_1, ..., \kappa_n)$ lies in the cone $\Gamma$. If we choose a family of suitable diffeomorphisms \ $ \varphi(t , . ): \mathbb{S}^n \to \mathbb{S}^n $ then
 $$
 X(t , x)= \rho(t, \varphi(t , x)) \varphi(t , x ),
 $$
 where \ $ \rho(t , . ) : \mathbb{S}^n \to \mathbb{R}^{+}$ is the radial function of \ $M_t$.  We have  
 $$
 \partial_t X  = \left( \langle\nabla\rho ,  \partial_t\varphi\rangle + \partial_t \rho \right) \varphi + \rho  \partial_t \varphi  $$
and the unit outer normal is given by 
$$\nu = {   \rho\varphi  - \nabla\rho \over \sqrt{|\nabla\rho|^2 + \rho^2}} \ .$$
 
 \medskip
 
\noindent  Using the fact that \ $ \partial_t \varphi $ is tangential to \ $S^n$ at \ $ \varphi$, it follows that 
 $$
 \left\langle  \partial_t X  , \nu \right\rangle = \left( \rho^2 + | \nabla \rho |^2  \right)^{-{1 \over 2}} \rho \partial_t \rho$$
 hence \ $\rho$ satisfies the initial value problem
 $$
\begin{cases} \displaystyle \partial_t\rho =  \mathcal{F}[\rho(t, .)]  \cr 
\rho(0 , x) = \rho_0(x) , \  x \in \mathbb{S}^n \end{cases} \eqno (2.10)
$$
where  the nonlinear operator $\mathcal{F}$ is defined on smooth functions  $\rho : \mathbb{S}^n \to (0, +\infty) $, such that  the matrix  $[a_{ij}]$ given in (2.4) lies in $M(\Gamma)$,  by 
$$ \mathcal{F}[\rho](x)=  \left( {1 \over F(a_{ij}(x))}  - f(\rho(x) x) \right) { \sqrt{ \rho^2(x) + | \nabla \rho(x) |^2  } \over \rho(x)} \ .   \eqno (2.11)$$
From now on, what we mean by admissible function   is a smooth function $\rho : [0, T]\times \mathbb{S}^n \to (0, + \infty)$ such that the matrix $[a_{ij}]$ defined by 
(2.4) lies in the cone $M(\Gamma)$ defined above.  Conversely, suppose that $\rho : [0, T]\times\mathbb{S}^n \to (0, +\infty)$ is an admissible  solution of (2.10) .   If we set 
$$X(t, x) = \rho(t, \varphi(t, x))\varphi(t, x) \ , \  (t, x) \in [0, T]\times  \mathbb{S}^n,$$
where $\varphi(t, .) : \mathbb{S}^n \to \mathbb{S}^n $ is a diffeomorphism satisfying the ODE 
$$\begin{cases} \partial_t\varphi(t, x) = Z(t, \varphi(t, x)  )
 \cr \varphi(0 , x) = x , \   x \in \mathbb{S}^n \end{cases} \eqno (2.12)$$
 with 
 $$Z(t, y) = -  \left( {1 \over F(a_{ij}(t,y))}  - f(\rho(t, y) y) \right) {\nabla \rho(t,y) \over \rho  \sqrt{|\nabla \rho(t,y)|^2 + \rho^2(t, y)}} ,  \ (t, y) \in  [0, T]\times\mathbb{S}^n,  \eqno (2.13) $$
 
 \medskip
 
\noindent  then  it is not difficult to see  that $X$  is a  starshaped embedding  which is  a solution of (1.2) with $X_0(x) = \rho_0(x)x. $

\medskip

The condition (2.6) implies that (2.10) is parabolic on admissible  functions $\rho$.  The classical theory of parabolic equations yields the existence and uniqueness  of a smooth admissible solution $\rho$ defined on a small  intervall $[0, T]$. From the classical theory of ordinary differential equations, there exists a family of diffeomorphisms $\varphi( t , .)$ defined on a small interval $[0, T]$ and satisfying (2.12). Thus by taking $X(t , x) = \rho(t , \varphi(t, x))\varphi(t, x)$ we obtain a solution of (1.2) defined on $[0, T]$.

\medskip

 Usually in order to get  high order estimates  it is useful to represent the hypersurface locally as graph over an open set $\Omega \subset \mathbb{R}^n$.  Locally, after rotating the coordinates axes , we may suppose that   $M$ is the graph of  a smooth function  $u : \Omega \to  \mathbb{R}$. Hence the metric of $M$,  the outer normal vector and the second fundamental form    can  be written respectively 
 
 $$g_{ij} =  \delta_{ij}+ D_iuD_ju \ , \   \    g^{ij} =  \delta_{ij}-{ D_iuD_ju \over 1 + |Du|^2} \eqno (2.14) $$
  $$\nu = {1 \over \sqrt{1+|Du|^2}} (Du, -1) \ ,  \eqno (2.15) $$
 $$h_{ij} = {D_{ij} u \over \sqrt{1+|Du|^2}}  \eqno (2.16) $$

 \medskip
 
 \noindent where $D_k $, $D_{ij}$ are the usual first  and second order derivatives in $\mathbb{R}^n$, and $Du = (D_1u , ... , D_nu)$. The principal curvatures  of $M$ are the eigenvalues of the symmetric matrix  $[a_{ij}]$ given by 
 $$ [a_{ij}] = [g^{ij} ]^{{1 \over 2}} [h_{ij} ] [g^{ij} ]^{{1 \over 2}} \eqno (2.17)$$
 where \ $ [ g^{ij} ]^{{1 \over 2}}$ is the positive square root of \ $ [ g^{ij} ]$. On ca compute  
 $$ a_{ij} = {1 \over v}\left( D_{ij}u - {D_iuD_luD_{jl}u \over v(1+v)} - { D_juD_luD_{il}u \over v(1+v)}  +  {D_iuD_juD_kuD_luD_{kl}u \over v^2(1+v)^2} \right)  \eqno (2.18)$$
 
 \noindent with $v= \sqrt{1 + |Du|^2} $. 
 
 \medskip

\noindent In this case equation (1.2) takes the forme 

$$ \partial_tu = - \left( { 1 \over F(a_{ij})} - f(x, u) \right) \sqrt{1 + |Du|^2}  \ .    \eqno (2.19)$$ 
 In what follows, what we mean by an admissible solution of (2.19) is a smooth function $u  :  [0, T]\times\Omega \to \mathbb{R}$ such that the matrix   $[a_{ij}]$ defined by (2.18) lies in the cone $M(\Gamma)$ defined above, and satisfying (2.19).

\section{ $C^1$-estimates and exponential decay}

\medskip
 
 In this section we prove $C^1$-estimates on  solutions  $\rho$ of  (2.10) and exponential decay of its derivatives $\partial_t\rho$.  First we prove  $C^0$-estimates.

\medskip

\newtheorem{prop}{Proposition}[section] 
\begin{prop} Suppose that $F$ satisfies conditions $(1.3)$-$(1.6)$ and  that $f$ satisfies conditions $(1.7)$-$(1.8)$. Let $\rho : [0, T]\times\mathbb{S}^n  \to (0, +\infty)$ be an admissible  solution of $(2.10)$. Then we have,  for all  $(t, x) \in [0, T]\times \mathbb{S}^n  $,
$$ R_1 \le  \rho(t, x)  \le R_2 \eqno (3.1) $$
where 
$$R_1 = \min\left(r_1, \min_{x \in \mathbb{S}^n}\rho_0(x) \right) \ \  \hbox{and} \  \    R_2 = \max\left(r_2, \max_{x \in \mathbb{S}^n}\rho_0(x) \right) $$
and where $r_1, r_2$  are as in (1.8). 
\end{prop}

\medskip

\begin{proof} Let $\rho : [0, T]\times\mathbb{S}^n  \to (0, +\infty)$ be an admissible  solution of (2.10). Let $(t_0 , x_0) \in [0, T] \times \mathbb{S}^n$ such that 
$$\rho(t_0 , x_0) = \max_{(t,x) \in [0, T]\times \in \mathbb{S}^n}\rho(t , x) . \eqno (3.2) $$

\noindent We want to prove 
$$ \rho(t_0 , x_0)  \le R_2.  \eqno (3.3) $$

\noindent If $t_0 = 0$, then 
$$\rho(t_0 , x_0) = \rho_0(x_0) \le R_2, $$
so (3.3) is proved in this case. Suppose now that $t_0 > 0.$  Then we have 
$$\partial_t\rho(t_0, x_0) \ge 0 \eqno (3.4) $$
 $$\nabla\rho(t_0, x_0)= 0   \eqno (3.5) $$

 \noindent and  the matrix 
  $$ \left[ \nabla_{ij}\rho(t_0, x_0) \right]   \ \text{is negative semi-definite}.   \eqno (3.6) $$
  
  \medskip
  
  \noindent It follows  from (3.5) and (3.6) that  the matrix $[a_{ij} ]$ defined by (2.4) satisfies 
in the sense of operators
$$ a_{ij}(t_0, x_0) \ge \rho^{-1}(t_0, x_0)\delta_{ij} . \eqno (3.7)$$

\medskip

\noindent Since  by (1.3) $F$ is monotone, then  by using (3.7)   we have at $(t_0, x_0)$
$$F(a_{ij}) \ge F(\rho^{-1}\delta_{ij}) = \rho^{-k}F(\delta_{ij})  =  \rho^{-k}, \eqno (3.8)$$
where we have used the fact that  $F$ is homogenous of degree $k$ and $F(\delta_{ij}) = 1$.  Using equation (2.10) and (3.8), we obtain 
$$\partial_t\rho(t_0, x_0)  \le  \rho^{k}(t_0, x_0) -  f(\rho(t_0, x_0) x_0).  \eqno (3.9)$$
Combining (3.4) and (3.9) gives 
$$f(\rho(t_0, x_0) x_0)  \le     \rho^{k}(t_0, x_0) . \eqno (3.10)$$

\noindent But from (1.7) and (1.8) we have  that if $ X\in \mathbb{R}^{n+1}$ satisfies $|X| > r_2$,  then  $f(X) > |X|^k$. So it follows from (3.10) that $ \rho(t_0, x_0) \le r_2 $.  This proves (3.3) since $r_2 \le R_2$. 

\medskip

It remains now to prove that  $\rho(t, x) \ge R_1$.  As before, if we let $(t_0, x_0)  \in [0, T]\times \mathbb{S}^n$ such that 
$$\rho(t_0 , x_0) = \min_{(t,x) \in [0, T]\times \in \mathbb{S}^n}\rho(t , x) ,$$
then in the same way as before, we prove that $\rho(t_0, x_0) \ge R_1$. This achieves the proof of Proposition 3.1.

\end{proof}

\medskip

We prove now the exponential decay of $\partial_t\rho$.  

\medskip

 \begin{prop} Assume  that $F$ satisfies conditions $(1.3)$-$(1.6)$ and  that $f$ satisfies conditions $(1.7)$-$(1.8)$. Let $\rho : [0, T]\times\mathbb{S}^n  \to (0, +\infty)$ be an admissible  solution of $(2.10)$.  We suppose  that $ \mathcal{F}[\rho_0] \ge 0$ if $k \le 1$ and $ \mathcal{F}[\rho_0] \le 0$ if $k > 1$, where the operator $\mathcal{F}$ is given by $(2.11)$,  and $k$ is the homogeneity degree of $F$.   Then we have, for any $(t, x ) \in [0, T]\times\mathbb{S}^n$, 
 $$\partial_t\rho(t, x) \ge 0 \  \  \hbox{if} \  \   k \le 1 $$
 and 
  $$\partial_t\rho(t, x) \le 0 \  \  \hbox{if} \  \   k > 1 . $$
  
 \noindent Moreover,  there exists  a positive constant $\lambda$ depending only on $f, r_1, r_2$ and $\rho_0$ such that,  for any $t \in [0, T]$,  we have 
 $$\max_{x \in \mathbb{S}^n} \left|\partial_t\rho(t, x)\right| \le{R_2 \over R_1} \max_{x \in \mathbb{S}^n} \bigl|\mathcal{F}[\rho_0](x)\bigr| e^{- \lambda t} \ , \eqno (3.11) $$
 where 
 $$R_1 = \min\left(r_1, \min_{x \in \mathbb{S}^n}\rho_0(x) \right) \ \  \hbox{and} \  \    R_2 = \max\left(r_2, \max_{x \in \mathbb{S}^n}\rho_0(x) \right) $$
and where $r_1, r_2$  are as in (1.8).

 \end{prop}
 
 \medskip

 \bigskip
 
 The proof of the above proposition is based on the following lemma which asserts that the function $ \rho^{-1}\partial_t\rho  $ satisfies a second order parabolic equation. 
 
 \medskip
 
 \newtheorem{lem}{Lemma}[section]
 
 \begin{lem}  Suppose that $F$ satisfies conditions $(1.3$)-$(1.6)$.  Let $\rho : [0, T]\times\mathbb{S}^n  \to  (0, +\infty)$ be an admissible  solution of  $(2.10)$ and set 
  $G =   \rho^{-1}\partial_t\rho$. 
 Then we have for some smooth functions $A_l, l= 1, ..., n$ $($ depending on $\rho $ and its derivatives $)$ , 
 $$\partial_t G=  \sum_{i, j=1}^nA_{ij} \nabla_{ij}G+  \sum_{l= 1}^nA_l\nabla_lG    -     { \sqrt{ \rho^2 + | \nabla \rho |^2 }  \over \rho^2} \left(  \rho \partial_{\rho}f  - f -  {k-1 \over F} \right)G $$
 where
   $$A_{ij}  = { 1 \over \rho^{2} F^2 } \sum_{l, m = 1}^n  \gamma_{il}F_{lm}gamma_{mj}  \eqno (3.12) $$
  and
  $$\gamma_{ij} = \delta_{ij} -{\nabla_i\rho\nabla_j\rho \over\sqrt{\rho^2 + |\nabla \rho|^2} \left(\rho + \sqrt{\rho^2 + |\nabla \rho|^2}\right)} \  . \eqno (3.13) $$
 
 \end{lem}

 \bigskip

 \begin{proof} We recall that  by  (2.10), $\rho$  satifies 
 
 $$\ \partial_t \rho =  \mathcal{F}[\rho]  \eqno (3.14) $$
 where
 $$ \mathcal{F}[\rho]  = \left( {1 \over F(a_{ij})}  - f(\rho x) \right) { \sqrt{ \rho^2 + | \nabla \rho |^2  } \over \rho }  \eqno (3.15)$$
and where $a_{ij}$ is given by (2.4). 

\medskip

\noindent In view of the definition of $G$ and (3.15) it will be usefull to work with the function $r = \log\rho$ instead of $\rho$. Equation (3.14) becomes then 
 $$\ \partial_t r =  \left( {1 \over F(a_{ij})}  - f(e^{r} x) \right) e^{-r} \sqrt{1 + | \nabla r|^2  }   \eqno (3.16) $$
where $a_{ij}$ takes the form

$$ a_{ij} = {e^{-r}b_{ij} \over \sqrt{1 + |\nabla r|^2}}  \eqno (3.17) $$
with
$$\begin{cases}  \displaystyle b_{ij} = \gamma_{il}(\delta_{lm} + \nabla_lr\nabla_mr -\nabla_{lm}r )\gamma_{mj} \cr 
\cr \displaystyle 
\gamma_{ij} = \delta_{ij} -{\nabla_ir\nabla_jr \over\sqrt{1 + |\nabla r|^2} \left(1 + \sqrt{1 + |\nabla r|^2}\right) } .  \end{cases} \eqno (3.18) $$

\bigskip

 Now,  we have 
 $$G= \rho^{-1}\partial_t\rho = \partial_tr =   \left( {1 \over F(a_{ij})}  - f(e^{r} x) \right)   e^{-r}\sqrt{1 + | \nabla r|^2  } , \eqno (3.19) $$
 so 
$$\partial_tG = - e^{-r}\sqrt{1 + | \nabla r|^2  }  \sum_{i, j=1}^n {F_{ij}\over F^2}\partial_ta_{ij} -  \sqrt{1 + | \nabla r|^2  }  \partial_{\rho}f(e^r x) \partial_tr $$
$$   +  \left( {1 \over F(a_{ij})}  - f(e^{r} x) \right) e^{-r} \left( -  \sqrt{1 + | \nabla r|^2  }\partial_t r   +   { \langle \nabla \partial_tr, \nabla r\rangle \over \sqrt{1 + | \nabla r|^2  } } \right)  .  \eqno(3.20) $$
   
\medskip

\noindent Using (3.17) and (3.18), one can check  that for some smooth functions $B_{ij}^l(t,x)$ ($l= 1, ..., n$), we have 
$$ \partial_ta_{ij} =- a_{ij}\partial_tr  -{ e^{-r} \over \sqrt{1 + | \nabla r|^2  } }\sum_{l,m=1}^n \gamma_{il}\gamma_{mj}\nabla_{lm}\partial_tr + \sum_{l= 1}^nB_{ij}^l \nabla_l\partial_tr ,  \eqno (3.21)$$

\noindent and since $\partial_tr = G$,  it follows from (3.20) and (3.21) that 

 $$\partial_tG = \sum_{i, j=1}^nA_{ij}\nabla_{ij}G +  \sum_{l= 1}^nA_l\nabla_lG  -  \partial_{\rho}f(e^r x) \sqrt{1 + | \nabla r|^2  }  G - G^2  $$
 $$+  \  e^{-r} \sqrt{1 + | \nabla r|^2  } \sum_{i,j=1}^n{F_{ij} \over F^2}a_{ij} G, \eqno (3.22) $$
 
 \medskip
 
 \noindent where
  $$A_{ij} ={ e^{-2r}  \over F^2 } \sum_{l, m=1}^n \gamma_{il}\gamma_{mj} F_{lm}$$
  and \  $A_l(t,x)$ ( $l=1, ..., n $) are smooth functions.   Since $F$ is homogeneous of degree $k$, then 
  
  $$ \sum_{i,j=1}^n{F_{ij} \over F^2}a_{ij}  = {k \over F},$$
  so it follows from (3.22) that 
  
  $$\partial_tG = \sum_{i, j=1}^nA_{ij}\nabla_{ij}G +  \sum_{l= 1}^nA_l\nabla_lG    -    \sqrt{1 + | \nabla r|^2  }  e^{-r} \left( e^r \partial_{\rho}f     -  {k \over F} \right) G - G^2 $$
  
 $$=  \sum_{i, j=1}^nA_{ij}\nabla_{ij}G +  \sum_{l= 1}^nA_l\nabla_lG    -      { \sqrt{\rho^2 + | \nabla \rho |^2 }  \over \rho^2} \left(  \rho\partial_{\rho}f  - f -  {k-1 \over F} \right)G. $$
 
 \noindent This achieves the proof Lemma 3.1. 
 
 \end{proof}

 \bigskip

 We need also the following lemma which is a well known version of the maximum principle for  parabolic equations. 
 
 \medskip
 
 \begin{lem} Let $G : [0, T]\times \mathbb{S}^n \to \mathbb{R}$ be a smooth function satisfying 
  $$\partial_t G\ge  \sum_{i, j=1}^nA_{ij}\nabla_{ij}G+  \sum_{l= 1}^nA_l\nabla_lG  + A G  \eqno (3.23)$$
  for some smooth functions $A, A_l, A_{ij}, ( l,i, j = 1, ..., n)$, such that the matrix $\left[ A_{ij} \right]$ is positive semi-definite.  Suppose 
  $$\min_{x \in \mathbb{S}^n}G(0, x) \ge 0, $$
  then  
  $$\min_{(t,x)\in [0, T]\times \mathbb{S}^n}G(t, x) \ge 0.$$
  \end{lem}
 
 \medskip
 
 \begin{proof} Let  $\lambda \in \mathbb{R}$ such that 
 $$\lambda < -   \max_{(t, x) \in [0, T]\times \mathbb{S}^n}|A(t , x)|  \  , \eqno (3.24) $$

 \noindent and consider the function $\widetilde{G}$ defined by $\widetilde{G}(t, x) = e^{\lambda t} G(t,x) $.    To prove the lemma it is  equivalent  to prove that 
   $$  \min_{(t, x) \in [0, T ]  \times\mathbb{S}^n}\widetilde{G}(t , x)  \ge 0. \eqno (3.25)$$

 \noindent By using (3.23), $\widetilde{G}$ satisfies 
   $$\partial_t \widetilde{G}\ge  \sum_{i, j=1}^nA_{ij}\nabla_{ij}\widetilde{G} +  \sum_{l= 1}^nA_l\nabla_l\widetilde{G}  + (\lambda + A )\widetilde{G} . \eqno (3.26)$$
 
   \noindent  Let  $(t_0, x_0) \in [0, T] \times \mathbb{S}^n$ such that 
 $$\widetilde{G}(t_0, x_0) = \min_{(t, x) \in [0, T] \times \mathbb{S}^n}\widetilde{G}(t , x) .$$
 
 \noindent We want to prove 
 $$\widetilde{G}(t_0, x_0) \ge 0.  \eqno (3.27)$$
 
 \noindent If $t_0 = 0$, then 
 
$$\widetilde{G}(t_0, x_0) = \widetilde{G}(0, x_0) = G(0, x_0)  \ge 0$$
and (3.27) is proved in this case. If $t_0 >0$, then 
$$\partial_t \widetilde{G}(t_0, x_0) \le 0 \eqno (3.28)$$
$$\nabla \widetilde{G}(t_0, x_0) =0  \eqno (3.29)$$
and the matrix 
$$\left[ \nabla_{ij}\widetilde{G}(t_0, x_0) \right]   \  \text{ is positive semi-definite}. \eqno (3.30) $$ 

\noindent It follows from (3.26), (3.28), (3.29) and (3.30) that 
$$ (\lambda + A(t_0, x_0))\widetilde{G}(t_0, x_0) \le 0$$
which implies that $\widetilde{G}(t_0, x_0) \ge 0$ since $ \lambda + A(t_0, x_0) < 0 $ by (3.24).  Thus (3.27) is proved and the lemma follows.

 \end{proof}
 
 \begin{proof}[Proof of Proposition 3.2] \    Let $G = \rho^{-1}\partial_t\rho$. Then by  Lemma 3.1 we have 
 
 $$ \partial_t G=  \sum_{i, j}A_{ij}\nabla_{ij}G +  \sum_{l= 1}^nA_l\nabla_lG  $$
 $$  -  \   { \sqrt{\rho^2 + | \nabla \rho |^2 }  \over \rho^2} \left(  \rho \partial_{\rho}f  - f -  {k-1 \over F} \right)G. \eqno (3.31)$$

 \noindent  By (1.3) (or equivalently (2.6)) the matrix $\left[F_{ij}\right]$ is positive definite. So it follows from (3.12)  that $[A_{ij}]$ is  positive semi-definite. We distinguish two cases :
 
 \medskip
 
 \noindent \underline{First case} :  $0 < k \le 1$  . Since $G$ satisfies (3.31) and $G(0, x)  =  \rho_0^{-1}(x)\partial_t\rho(0, x) =  \rho_0^{-1}(x) \mathcal{F}[\rho_0](x) \ge 0$ by hypothesis,  then by Lemma 3.2 we have for any $t \in [0, T]$, 
 $$\min_{x \in \mathbb{S}^n} G(t, x) \ge 0.  \eqno (3.32) $$
 In particular, (3.32) implies  that $\partial_t\rho \ge 0$ since $\partial_t\rho = \rho G$.   Now we have,  since $\rho$ satisfies (2.10), 
$$ G = \rho^{-1}\partial_t\rho = \left({1\over F(a_{ij})} - f(\rho x)\right){\sqrt{\rho^2 + |\nabla \rho|^2 } \over \rho^2}\ , $$
 so  it follows   from (3.32) that 
 $${1 \over F(a_{ij})} \ge  f(\rho x) $$
 which implies that the last term in (3.31)  is  bounded from below as 
  $$  { \sqrt{\rho^2 + | \nabla \rho |^2 }  \over \rho^2}  \left(  \rho \partial_{\rho}f  - f -  {k-1 \over F} \right) \ge
      { \sqrt{\rho^2 + | \nabla \rho |^2 }  \over \rho^2}  \left(  \rho \partial_{\rho}f - kf  \right) .  \eqno (3.33) $$
   
   \noindent  Since $f$ satisfies (1.7), then $ \rho \partial_{\rho}f - kf   >0$, and since $R_1 \le \rho(t , x) \le R_2 $ by  Proposition 3.1,    we deduce  that 
   $$ \rho \partial_{\rho}f - kf   \ge \delta_0   \eqno (3.34) $$
      for some constant $\delta_0 >0$  depending  only on   $f, R_1$ and $R_2$.  It  follows from (3.33)  and (3.34) by using Proposition 3.1 that 
   $$ { \sqrt{\rho^2 + | \nabla \rho |^2 }  \over \rho^2}  \left(  \rho \partial_{\rho}f  - f -  {k-1 \over F} \right) \ge  {\delta_0 \over R_2} \ .  \eqno (3.35) $$
   
   \noindent By setting $\lambda =  {\delta_0 \over R_2}  $  and $\widetilde{G}(t, x) = e^{\lambda t}G(t,x)$, it follows from (3.31) that $\widetilde{G}$ satisfies 
    $$\partial_t \widetilde{G}=  \sum_{i, j=1}^nA_{ij}\nabla_{ij}\widetilde{G}+  \sum_{l= 1}^nA_l\nabla_l\widetilde{G}  $$    
 $$   - \   { \sqrt{\rho^2 + | \nabla \rho |^2 }  \over \rho^2}   \left(  \rho \partial_{\rho}f  - f -  {k-1 \over F} \right) \widetilde{G}  +   \lambda \widetilde{G}  $$
 which gives  by using (3.35) and the fact that $\widetilde{G} \ge 0$, 
  $$\partial_t \widetilde{G}\le   \sum_{i, j=1}^nA_{ij}\nabla_{ij}\widetilde{G}+  \sum_{l= 1}^nA_l\nabla_l\widetilde{G}.    \eqno (3.36) $$  
It follows from (3.36)  by applying Lemma 3.2 to the function $-  \widetilde{G} + \displaystyle  \max_{x \in \mathbb{S}^n}\widetilde{G}(0,x)$ that 
  $$ -  \widetilde{G} + \displaystyle  \max_{x \in \mathbb{S}^n}\widetilde{G}(0,x) \ge 0$$
    which implies 
  $$ \max_{x \in \mathbb{S}^n}G(t,x) \le e^{-\lambda t} \max_{x \in \mathbb{S}^n}G(0,x) . \eqno (3.37) $$
  
 \noindent But from the definition of $G$ we have 
  $$\partial_t\rho =\rho G,  \eqno (3.38) $$ 
  so it follows from (3.37) and (3.38) since $\partial_t\rho \ge 0$ and  $R_1 \le \rho \le R_2$ by Proposition 3.1,  that 
  $$|\partial_t\rho | \le R_2 e^{- \lambda t} \max_{x \in \mathbb{S}^n}G(0,x) = R_2 e^{- \lambda t} \max_{x \in \mathbb{S}^n}\left({\mathcal{F}[\rho_0](x) \over \rho_0(x)}\right) \le 
   {R_2 \over R_1} e^{- \lambda t} \max_{x \in \mathbb{S}^n} \mathcal{F}[\rho_0](x)  . $$
  This proves Proposition 3.2  in the case $ 0 < k \le 1$.

 \medskip

 \noindent \underline{Second case} :  $ k > 1$.  \  Since $G$ satisfies (3.31) and $G(0, x)  =  \rho_0^{-1}(x)\partial_t\rho(0, x) =  \rho_0^{-1}(x) \mathcal{F}[\rho_0](x) \le 0$ by hypothesis,  then by Lemma 3.2 we have for any $t \in [0, T]$, 
 $$\max_{x \in \mathbb{S}^n} G(t, x) \le 0.  \eqno (3.39) $$
 In particular, (3.39) implies  that $\partial_t\rho \le 0$ since $\partial_t\rho = \rho G$.   Now we have,  since $\rho$ satisfies (2.10), 
$$ G = \rho^{-1}\partial_t\rho = \left({1\over F(a_{ij})} - f(\rho x)\right){\sqrt{\rho^2 + |\nabla \rho|^2 } \over \rho^2}\ , $$
 so  it follows   from (3.39) that 
 $${1 \over F(a_{ij})} \le  f(\rho x) $$
 which implies that the last term in (3.31)  is  bounded from below as 
  $$  { \sqrt{\rho^2 + | \nabla \rho |^2 }  \over \rho^2} \left(  \rho \partial_{\rho}f  - f -  {k-1 \over F} \right) \ge
      { \sqrt{\rho^2 + | \nabla \rho |^2 }  \over \rho^2}  \left(  \rho \partial_{\rho}f - kf  \right) .  \eqno (3.40) $$
   
   \noindent  Since $f$ satisfies (1.7), then $ \rho \partial_{\rho}f - kf   >0$, and since $R_1 \le \rho(t , x) \le R_2 $ by  Proposition 3.1,    we deduce  that 
   $$ \rho \partial_{\rho}f - kf   \ge \delta_0   \eqno (3.41) $$
      for some constant $\delta_0 >0$  depending  only on   $f, R_1$ and $R_2$.  It  follows from (3.40)  and (3.41) by using Proposition 3.1 that 
   $$ { \sqrt{\rho^2 + | \nabla \rho |^2 }  \over \rho^2}  \left(  \rho \partial_{\rho}f  - f -  {k-1 \over F} \right) \ge  {\delta_0 \over R_1} \ .  \eqno (3.42) $$
   
   \noindent By setting $\lambda =  {\delta_0 \over R_1}  $  and $\widetilde{G}(t, x) = e^{\lambda t}G(t,x)$, it follows from (3.31) that 
    $$\partial_t \widetilde{G}=  \sum_{i, j=1}^nA_{ij}\nabla_{ij}\widetilde{G}+  \sum_{l= 1}^nA_l\nabla_l\widetilde{G}  $$    
 $$   - \  { \sqrt{\rho^2 + | \nabla \rho |^2 }  \over \rho^2}  \left(  \rho \partial_{\rho}f  - f -  {k-1 \over F} \right) \widetilde{G}  +   \lambda \widetilde{G}  $$
 which gives  by using (3.42) and the fact that $\tilde{G} \le 0$, 
  $$\partial_t \widetilde{G}\ge   \sum_{i, j=1}^nA_{ij}(t, x)\nabla_{ij}\widetilde{G}+  \sum_{l= 1}^nA_l(t, x)\nabla_l\widetilde{G}.    \eqno (3.43) $$  
It follows from (3.43)  by applying Lemma 3.2 to the function $ \widetilde{G} - \displaystyle  \min_{x \in \mathbb{S}^n}\widetilde{G}(0,x)$ that 
  $$ \widetilde{G} - \displaystyle  \min_{x \in \mathbb{S}^n}\widetilde{G}(0,x) \ge 0$$
    which implies 
  $$ \min_{x \in \mathbb{S}^n}G(t,x) \ge e^{-\lambda t} \min_{x \in \mathbb{S}^n}G(0,x) . \eqno (3.44) $$
  
 \noindent But from the definition of $G$ we have 
  $$\partial_t\rho =\rho G,  \eqno (3.45) $$ 
  so it follows from (3.44) and (3.45) since $\partial_t\rho \le 0$ and  $\rho \le R_2$ by Proposition 3.1,  that 
  $$|\partial_t\rho | \le - R_2 e^{- \lambda t} \min_{x \in \mathbb{S}^n}G(0,x) =   R_2 e^{- \lambda t} \max_{x \in \mathbb{S}^n}|G(0,x)| =  R_2 e^{- \lambda t} \max_{x \in \mathbb{S}^n}\left({\bigl|\mathcal{F}[\rho_0](x)\bigr| \over \rho_0(x)}\right)$$
  $$ \le 
   {R_2 \over R_1} e^{- \lambda t} \max_{x \in \mathbb{S}^n} \bigl|\mathcal{F}[\rho_0](x)\bigr|  . $$
\noindent The proof of Proposition 3.2 is then complete.

 \end{proof}

 Now we are in position to prove $C^1$-estimates on the function $\rho$. 
 
 \medskip
 
  \begin{prop} Supoose that $F$ satisfies conditions $(1.3)$-$(1.6)$  and that $f$ satisfies conditions $(1.7)$-$(1.8)$. Let $\rho : [0, T]\times\mathbb{S}^n  \to \mathbb{R}^{+}$ be an admissible   solution of $(2.10)$.  We suppose  that $ \mathcal{F}[\rho_0] \ge 0$ if $k \le 1$ and $ \mathcal{F}[\rho_0] \le 0$ if $k > 1$, where the operator $\mathcal{F}$ is given by $(2.11)$,  and $k$ is the homogeneity degree of $F$. Then there exists a positive constant  $C$ depending only on $f, r_1, r_2$ and $\rho_0$ such that 
  $$\max_{(t,x) \in [0, T] \times\mathbb{S}^n} |\nabla\rho(t, x)| \le C ,$$
  where $r_1$ and $r_2$ are as in (1.8). 
  \end{prop}
  
  \medskip
  
 \begin{proof} 
 As in the proof of Lemma 3.1,  we introduce the function $r =  \log\rho$. We have then 
  $$\ \partial_t r =  \left( {1 \over F(a_{ij})}  - f(e^{r} x) \right) e^{-r}\sqrt{ 1 + | \nabla r|^2   }   \eqno (3.46) $$
where we recall that  $a_{ij}$ takes the form

$$ a_{ij} = {e^{-r} b_{ij}  \over \sqrt{1 + |\nabla r|^2} }  \eqno (3.47) $$
with
$$\begin{cases}   \displaystyle  b_{ij} = \gamma_{il}(\delta_{lm} + \nabla_lr\nabla_mr -\nabla_{lm}r )\gamma_{mj} \cr\cr

\displaystyle \gamma_{ij} = \delta_{ij} -{\nabla_ir\nabla_jr \over\sqrt{1 + |\nabla r|^2} \left(1 + \sqrt{1 + |\nabla r|^2}\right) } .  \end{cases} \eqno (3.48)$$

\bigskip
 
 Set $H = {1 \over 2} |\nabla r|^2$, and  let $(t_0, x_0) \in [0, T]\times\mathbb{S}^n$ such that 
 $$H(t_0, x_0) =  \max_{(t,x) \in [0, T]\times\mathbb{S}^n} H(t , x) .$$ 
 Let $\{ e_1, ..., e_n \}$  be an orthonormal frame in a neighborhood  of $x_0$ such that  $\nabla_i(e_j) = 0 $ at $x_0$, for $i,j = 1, ..., n$. 
 
 \medskip

 \noindent If $t_0 = 0$, then 
$$ H(t_0, x_0) = H(0, x_0)  = \max_{x \in \mathbb{S}^n} H(0 , x) . \eqno (3.49) $$

\noindent If $t_0 > 0$, then 
$$\partial_tH(t_0, x_0) \ge 0 \eqno (3.50)$$

$$\nabla_i H(t_0, x_0) = 0, \ i=1, ..., n  \eqno (3.51) $$
and the matrix 
$$ \left[ \nabla_{ij}H(t_0, x_0) \right]  \  \text{is negative semi-definite. }  \eqno (3.52)$$

\medskip

 In what follows, to simplify the notataion we shall write $F$ instead of $ F(a_{ij})$, and $f$ instead $f(e^r x)$.  We have  at $(t_0, x_0) $, by using (3.51), 
$$\partial_t H = \langle \nabla\partial_t r, \nabla r \rangle = \left\langle \nabla\left( \left({1 \over F }  - f  \right)e^{-r} \sqrt{1 + | \nabla r|^2 } \right) ,   \nabla r  \right\rangle  $$
$$ =  - e^{-r}\sqrt{ 1 + | \nabla r|^2   }  \sum_{i, j=1}^n {F_{ij}\over F^2}\langle\nabla a_{ij} , \nabla r \rangle -  2 \sqrt{ 1 + | \nabla r|^2} \partial_{\rho} f H  $$
$$  -  \      \sqrt{ 1 + | \nabla r|^2  }\langle \nabla f , \nabla r\rangle   - 2   \left( {1 \over F }  - f  \right)e^{-r}\sqrt{ 1 + | \nabla r|^2  }  H   \  .  \eqno (3.53)$$
 
 \medskip
 
 \noindent \noindent Using (3.47) and (3.48), one can check  that for some smooth functions $B_{ij}^l(t,x)$ ($l= 1, ..., n$), we have, for any $\alpha= 1, ..., n$, at $ (t_0, x_0)$, 
$$ \nabla_{\alpha}a_{ij} = - {e^{-r} \over \sqrt{1 + | \nabla r |^2 }}\sum_{l,m =1}^n \gamma_{il}\gamma_{mj}\nabla_{\alpha lm}r + \sum_{l= 1}^nB_{ij}^l \nabla_{\alpha l}r - a_{ij}\nabla_{\alpha}r . $$

\noindent It follows that, at  $(t_0, x_0)$, 

 $$ \langle\nabla a_{ij} , \nabla r \rangle =  \sum_{\alpha=1}^n   \nabla_{\alpha}a_{ij}\nabla_{\alpha}r$$
 $$ =  
  -  {e^{-r}\over \sqrt{1 + | \nabla r |^2 }} \sum_{\alpha, l,m =1}^n \gamma_{il}\gamma_{mj}\nabla_{\alpha lm}r \nabla_{\alpha}r 
  - 2a_{ij}H.   \eqno (3.54)$$
   The formula for commuting the order of covariant differentiation  gives  at $(t_0, x_0)$ 
$$ \nabla_{\alpha lm}r  =   \nabla_{lm\alpha}r + \delta_{\alpha m}\nabla_{l}r - \delta_{lm}\nabla_{\alpha}r . \eqno (3.55)$$
 
 \noindent Combining (3.54) and (3.55) we get at $(t_0, x_0)$ 
 $$  \langle\nabla a_{ij} , \nabla r \rangle =  - {e^{-r} \over \sqrt{1 + | \nabla r |^2 }} \sum_{\alpha, l,m =1}^n \gamma_{il}\gamma_{mj}\nabla_{lm\alpha}r\nabla_{\alpha}r  $$
$$  -  \  {e^{-r}\over \sqrt{1 + | \nabla r |^2 }} \sum_{ l,m =1}^n \gamma_{il}\gamma_{mj} \nabla_{l}r\nabla_m r $$
$$ + \  2 {e^{-r}\over \sqrt{1 + | \nabla r |^2 }} \sum_{ l =1}^n \gamma_{il}\gamma_{lj} H - 2a_{ij}H.  \eqno (3.56)$$

  \noindent But we have at $(t_0, x_0)$ 
  $$\nabla_{lm}H ={1\over 2} \nabla_{lm}\left(|\nabla r|^2\right) =  \sum_{\alpha =1}^n \nabla_{lm\alpha}r \nabla_{\alpha}r  +   \sum_{\alpha  =1}^n \nabla_{l \alpha }r \nabla_{m \alpha }r .    \eqno (3.57)$$
  Hence it follows from (3.50), (3.53),  (3.56)  and (3.57) that, at $(t_0, x_0)$, 
  
  $$ 0 \le e^{-2r} \sum_{ i , j =1}^nA_{ij}\nabla_{ij}H  - e^{-2r} \sum_{\alpha,  l , m =1}^nA_{lm}\nabla_{l \alpha }r \nabla_{m \alpha }r $$
$$  +  \   2 e^{-r} \sqrt{ 1 + | \nabla r|^2   }   \sum_{i, j = 1}^n{F_{ij}\over F^2}a_{ij} H   + e^{-2r}  \sum_{i , j = 1}^nA_{ij}\nabla_ir\nabla_jr  -  2 e^{-2r} \text{Trace}\left[A_{ij}\right] H$$
$$  - 2 \sqrt{ 1 + | \nabla r|^2  } \partial_{\rho} f H   -       \sqrt{1 + | \nabla r|^2} \langle \nabla f , \nabla r\rangle -  \  2   \left( {1 \over F }  - f  \right)e^{-r}\sqrt{ 1 + | \nabla r|^2 }  H ,   \eqno (3.58) $$

 where 
 $$A_{ij} = \sum_{l, m = 1}^n{F_{lm}\over F^2}\gamma_{il}\gamma_{mj} .$$
 Since $[F_{ij}]$  is positive definite, then  $[ A_{ij}]$  is positive semi-definite. So we have  at $(t_0, x_0)$,  by using (3.52), 
 $$  \sum_{ i , j =1}^nA_{ij}\nabla_{ij}H \le 0, \eqno (3.59) $$ 
 $$     \sum_{\alpha,  l , m =1}^nA_{lm}\nabla_{l \alpha }r \nabla_{m \alpha }r \ge 0 \eqno (3.60)$$
 and
 $$\sum_{i , j = 1}^nA_{ij}\nabla_ir\nabla_jr  -  2\  \text{Trace}\left[A_{ij}\right] H \le 0.  \eqno (3.61) $$
 
 Since $F$ is homogenous of degree $k$,  we have also 
 $$\sum_{i, j = 1}^n{F_{ij}\over F^2}a_{ij} = {k \over F} . \eqno (3.62)$$
 Thus we get from (3.58),  (3.59), (3.60), (3.61), (3.61) and (3.62), at $(t_0, x_0)$ 
 
$$ 0  \le     \ 2 e^{-r} \sqrt{1 + | \nabla r|^2   } {k \over F}  H -   2 \sqrt{1 + | \nabla r|^2 } \partial_{\rho} f H  $$
$$   -  \  2   \left( {1 \over F }  - f  \right)e^{-r} \sqrt{1 + | \nabla r|^2   }  H  -      \sqrt{ 1 + | \nabla r|^2}\langle \nabla f , \nabla r\rangle  . \eqno (3.63) $$
But by Proposition 3.2 we have $\partial_t\rho \ge 0$ if $k\le 1$, and $\partial_t\rho \le 0$ if $k > 1$. This implies,  since $\rho$ satifies (2.10), that  
$\displaystyle {1 \over F(a_{ij})} - f(\rho x) \ge 0 $ if $k \le 1$, and $\displaystyle {1 \over F(a_{ij})} - f(\rho x) \le 0 $ if $k > 1$. That is, 
$${k-1 \over F(a_{ij})} \le (k-1) f(\rho x)$$
Hence it follows  from (3.63) that  at $(t_0, x_0)$ 
 $$   2 \left( e^{r} \partial_{\rho} f   - kf \right) H  \le   e^{r} \langle \nabla f , \nabla r\rangle . \eqno (3.64) $$
 
 \medskip
 
 \noindent By (1.7) we have  $\rho \partial_{\rho} f(\rho x)   - kf(\rho x ) > 0$, which implies that 
 $$ \delta_0 = \min_{(\rho , x) \in [ R_1  ,  R_2]\times\mathbb{S}^n} \bigl( \rho\partial_{\rho}f(\rho x)- kf(\rho x )  \bigr) >  0,$$
 where $R_1$ and $R_2$ are defined in Proposition 3.1. Since $ R_1 \le \rho(t, x) \le R_2$ by Proposition 3.1, then   $ e^{r} \partial_{\rho} f   - kf  \ge \delta_0$. Thus it follows from (3.64)  at $(t_0, x_0)$ 
$$2 \delta_0 H \le  e^r \langle \nabla f , \nabla r\rangle \le  R_2 |\nabla f| |\nabla r|  = R_2 |\nabla f| \sqrt{2}\sqrt{H} $$
that is 
$$H(t_0, x_0) \le {C_0^2R_2^2 \over 2 \delta_0^2} \ ,  \eqno (3.66)$$ where 
 $$ C_0 =     \sup_{R_1 \le |y| \le R_2 }|\nabla f(y)| .$$
 
 \noindent It follows from (3.49) and (3.66) that 
$$ H(t_0, x_0)  \le  \max\left(   \max_{x \in \mathbb{S}^n} H(0 , x) , {C_0^2R_2^2 \over 2 \delta_0^2}\right) .$$

 \noindent This ends the proof of Proposition 3.3.
 
 \end{proof}

\section{$C^2$-estimates and proof of the main results}

To get $C^2$-estimates we need to controll the principal curvatures.

  \begin{prop} Suppose that $F$ satisfies conditions $(1.3)$-$(1.6)$ and  that $f$ satisfies conditions $(1.7)$-$(1.8)$. Let $\rho : [0, T]\times\mathbb{S}^n  \to (0, +\infty)$ be an admissible  solution of (2.10).   We suppose  that 
  $$  \begin{cases} \displaystyle \mathcal{F}[\rho_0] \ge 0 \   \    \hbox{if } \  \   k \le 1 \cr\cr 
   \displaystyle 0 \le -  \mathcal{F}[\rho_0] \le  {k R_1 \over (k+1)R_2}\min_{R_1 \le |Y| \le R_2} f(Y) \  \ \hbox{if } \ \  k > 1, \end{cases} \eqno (4.1) $$
   \medskip
  \noindent  where the operator $\mathcal{F}$ is given by $(2.11)$,   $k$ is the homogeneity degree of $F$, and  
 $$R_1 = \min\left(r_1, \min_{x \in \mathbb{S}^n}\rho_0(x) \right) \ \  ,  \  \    R_2 = \max\left(r_2, \max_{x \in \mathbb{S}^n}\rho_0(x) \right) $$
with  $r_1, r_2$   as in (1.8). Then there exists a positive constant  $C$ depending only on $f , r_1, r_2 $ and $\rho_0$ sucht that 
  $$\max_{(t,x) \in [0, T] \times\mathbb{S}^n} \max_{1\le i \le n}|\kappa_i(t, x)| \le C ,$$
  where $\kappa_1, ... , \kappa_n $  are the principal curvatures of  the hypersurface $M_t$  parametrized by $X(t ,  x) =  \rho(t , x) x$. 
  
  \end{prop}

\medskip

\begin{proof} Define the function $h : [0, T]\times \mathbb{S}^{n} \to \mathbb{R}$ by 
$$ h(t , x) = \log { {\displaystyle \max_{1 \le i \le n}\kappa_i(t ,x)}  \over \langle X(t ,x), \nu(t , x)\rangle} \eqno (4.2) $$
where  $\kappa_1, ... , \kappa_n$  are the principal curvatures of  the hypersurface  $M_t$ parametrized by $X(t, x) = \rho(t, x)x$, and $\nu(t,.)$ is its outer normal vector.  First we shall give   an upper bound on the function $h$.   Let  $(t_0 , x_0 ) \in  [0, T]\times\mathbb{S}^n$ the point where $h $ achieves its maximum on  $ [0, T]\times\mathbb{S}^n$, that is, 
$$h(t_0 ,  x_0) =   \max_{(t,x)  \in [0, T]\times\mathbb{S}^n}h(t , x) =   \max_{(t,x)  \in [0, T]\times \mathbb{S}^n}  \log { {\displaystyle \max_{1 \le i \le n}\kappa_i(t ,x)}  \over \langle X(t , x), \nu(t , x)\rangle}.$$

\noindent We want to prove that 
$$h(t_0, x_0) \le C_0 , \eqno (4.3) $$
\medskip

\noindent where the constant $C_0$ depends only on $f, r_1, r_2$ and $\rho_0$.   If $t_0 = 0$, then $h(t_0, x_0) = h(0, x_0) $, 
and (4.3) is trivially satisfied in this case. From now on, we suppose that $t_0 >0$. 
Without loss of generality, we may suppose that $x_0$ is the south pole of $\mathbb{S}^n$.  Let $\Sigma$ the tangent hyperplane to $M_{t_0}$ at the point $Z_0 = X(t_0, x_0)$. Then near $(t_0, Z_0)$, the family of hypersurfaces $M_t$ can be represented as  the graph of a smooth function $u$ defined on a neighbourhood of $(t_0, Z_0)$ in $[0, T] \times \Sigma$. Thus the function $u$ is  an admissible solution of (2.18) (see section 2). 
 
 By choosing a new coordinate system in the hyperplane $\Sigma$, with origin at the point $Z_0$, then  in the coordinate parallel to the new ones with centre at the original origin,  denoted by $x_1, ..., x_n$, we have 
 $$Z_0 = (a_1, ... , a_n, -a) \ , \  \text{for some constants } \ a_1, ..., a_n , a, \  \ \text{with} \  a >0, $$
 and 
 $$X(t, x) = (a_1, ... , a_n, -a)+ (x , u(t , x))  \   \    \text{with} \  u(t_0 , 0) = 0 .$$
 
 \medskip
 
\noindent  By formula (2.16) of section 2, we have 
 $$ \nu = {1 \over v}(Du, -1)   \eqno (4.4)$$
 and
 $$  \langle X ,  \nu \rangle = {1 \over v} \left( a - u +  \sum_{k=1}^n(x_k + a_k)D_ku\right)   \ , \eqno (4.5)$$
where 
$$v= (1 + |Du|^2)^{1/2} . \eqno (4.6)$$

\noindent By our choice of coordinates we have  

$$u(t_0, 0) = 0 \eqno (4.7)$$
and 
$$Du(t_0, 0) = (0, ..., 0)  .  \eqno (4.8) $$

\medskip

\noindent By rotating the new $x_1, ... , x_n$ coordinates, we may suppose that  $\displaystyle \max_{1 \le i \le n}\kappa_i(t_0 , x_0)$ occurs in the $x_1$-direction. We have then by using formula (2.17)  and (4.8) 
$$\max_{1 \le i \le n}\kappa_i(t_0 , x_0)   = \kappa_1(t_0 , x_0)  = {D_{11} u(t_0, 0) \over v(t_0, 0)(1 + \left( D_1 u(t_0, 0)\right)^2)} $$
$$ = D_{11}u(t_0, 0).$$

On a neighborhood of $(t_0, 0)$ define the function $H$ by 
$$H  = \log\left( {D_{11}u \over \varphi v (1 + (D_1u)^2)} \right)  $$
where 
$$\varphi = \langle X , \nu \rangle  =  {1 \over v} \left( a - u +  \sum_{k=1}^n(x_k + a_k)D_ku\right) $$
Thus we have 
$$H(t_0 , 0) = h(t_0, x_0) = \max_{(t,x) \in [0 , T]\times\mathbb{S}^n}h(t , x).   \eqno (4.9)$$
We  will give an upper bound on $H(t_0, 0)$.  By our choice of coordinates we have 
$$ D_{1\alpha}u(t_0 , 0) = 0 \  \text{for} \ \alpha > 1 ,  \eqno (4.10)$$
so by rotating the $x_2, ..., x_n$ coordinates, we may suppose that the matrix $D^2u(t_0, 0)$ is diagonal and that $D_{11}u(t_0 , 0) > 0 $.

\noindent We have,  since $H$ attains a local   maximum at $(t_0 , 0)$, that 
$$DH(t_0, 0) = 0 \eqno (4.11)  $$
and 
$$\partial_tH(t_0, 0) \ge 0 \eqno (4.12)$$
since $t_0 > 0$.  On the other hand, we have 
$$D_{\alpha}H = {D_{11\alpha}u \over D_{11}u} -  {D_{\alpha}v \over v} - {2D_1uD_{1\alpha}u \over 1 + (D_{1}u)^2} - {D_{\alpha}\varphi \over \varphi}$$
and
$$D_{\alpha}\varphi =\sum_{k=1}^n {(a_k + x_k)D_{\alpha k}u \over v} - {\varphi D_{\alpha}v \over v} .$$

\noindent But by using (4.8)  and (4.10) we have  at $(t_0 , 0)$ 
$$D_{\alpha} v = \sum_{k=1}^n{D_kuD_{\alpha k}u \over v} = 0,$$
so
$$D_{\alpha}\varphi   =  a_{\alpha}D_{\alpha\alpha}u  $$
and 
$$ D_{\alpha}H =  {D_{11\alpha}u \over D_{11}u} -   {a_{\alpha}D_{\alpha\alpha}u \over \varphi}  $$
which together with (4.11) give at $(t_0, 0)$, 
$${D_{11\alpha}u \over D_{11}u} -   {a_{\alpha}D_{\alpha\alpha}u \over \varphi}= 0. \eqno (4.13) $$

\noindent Differentiating once again, we get  at $(t_0 , 0)$
$$ D_{\alpha\alpha}v=  \left(D_{\alpha\alpha}u\right)^2$$
and 
$$D_{\alpha}\left({D_{\alpha}\varphi \over \varphi}\right) = {1 \over \varphi}\left( D_{\alpha\alpha}u + \sum_{k= 1}^n a_kD_{\alpha\alpha k}u\right) -  
{(a_{\alpha}D_{\alpha\alpha}u)^2 \over \varphi^2} - \left(D_{\alpha\alpha}u\right)^2.$$
So
$$D_{\alpha\alpha}H = {D_{11\alpha\alpha}u \over D_{11}u} - \left({D_{11\alpha}u \over D_{11}u}\right)^2 - 2 \left(D_{1\alpha}u \right)^2  + {(a_{\alpha}D_{\alpha\alpha}u)^2 \over \varphi^2}$$
$$- {1 \over \varphi}\left( D_{\alpha\alpha}u + \sum_{k= 1}^n a_kD_{\alpha\alpha k}u\right) $$
at $(t_0, 0)$.  And using (4.13) we obtain then 
$$D_{\alpha\alpha}H = {D_{11\alpha\alpha}u \over D_{11}u}  - 2 \left(D_{1\alpha}u \right)^2 - {1 \over a}\left( D_{\alpha\alpha}u + \sum_{k= 1}^n a_kD_{\alpha\alpha k}u\right) \eqno (4.14)
$$ at $(t_0, 0)$ for $\alpha = 1 , ..., n$, where de have used  the fact  that $a = \varphi(t_0, 0)$.

\medskip

Now if we differentiate equation (2.19) in the $x_1$ direction, we get 
$$D_1\partial_t u = -{ 1 \over \sqrt{1+ |Du|^2}}\left( {1 \over F} - f\right)\sum_{k=1}^nD_kuD_{k1}u $$
$$ + \  { \sqrt{1+ |Du|^2} \over F^2}\sum_{i, j =1}^n F_{ij}D_1a_{ij} \  + \   \sqrt{1+ |Du|^2 } \left(D_1f + D_{n+1}f D_1u \right) .$$
Differentiating once again in the $x_1$ direction and using (4.7), (4.8) and (4.10) we get at $(t_0, 0)$
$$D_{11}\partial_t u = -\left( {1 \over F} - f\right)\left(D_{11}u \right)^2 + { 1  \over F^2}\sum_{i, j =1}^n F_{ij}D_{11}a_{ij}  - {2 \over F^3}\left(\sum_{i,j=1}^nF_{ij}D_1a_{ij}\right)^2
$$
$$ + \   {1 \over F^2}\sum_{j,j,r,s=1}^nF_{ij,rs}D_1a_{ij}D_1a_{rs} + D_{11}f  + D_{n+1}f D_{11}u . \eqno (4.15) $$
But since $\log F$ is concave, we have
$$ - {2 \over F^3}\left(\sum_{i,j=1}^nF_{ij}D_1a_{ij}\right)^2 +    {1 \over F^2}\sum_{j,j,r,s=1}^nF_{ij,rs}D_1a_{ij}D_1a_{rs}  \le 0 , $$
so  it follows from (4.15) that 
$$D_{11}\partial_t u \le  -\left( {1 \over F} - f\right)\left(D_{11}u \right)^2 + { 1  \over F^2}\sum_{i, j =1}^n F_{ij}D_{11}a_{ij} +   D_{11}f  + D_{n+1}f D_{11}u   \eqno (4.16) $$ 
at $(t_0, 0)$. Now from the definition of the matrix $[a_{ij}]$ in (2.17), we have at $(t_0, 0)$ by using (4.7) and (4.8), 
$$D_{11}a_{ij}  = D_{11ij}u-  (D_{11}u)^2D_{ij}u -2 D_{1i}uD_{1j}uD_{11}u, $$
and since  $D^2u$ is  diagonal at $(t_0 , 0)$, then we have at this point 
$$D_{11}a_{11}  = D_{1111}u - 3(D_{11}u)^3 \eqno (4.17) $$
and 
$$D_{11}a_{\alpha\alpha}  = D_{11\alpha\alpha}u - D_{\alpha\alpha}u(D_{11}u)^2 \eqno (4.18) $$
for $\alpha = 2, ..., n$.  Combining (4.16), (4.17) and (4.18) we obtain,  since $\left[F_{ij}\right]$ is diagonal at $(t_0, 0)$, 
$$ D_{11}\partial_tu  \le - \left({1 \over F} - f\right)  (D_{11}u)^2 + {1 \over F^2} \left( \sum_{\alpha=1}^n F_{\alpha\alpha}D_{11\alpha\alpha}u  - (D_{11}u)^2\sum_{\alpha = 2}^n F_{\alpha\alpha}D_{\alpha\alpha}u  \right)$$
$$ -  \   3 {F_{11} \over F^2}(D_{11}u)^3 + D_{11}f  + D_{n+1}f D_{11}u . \eqno (4.19)$$
But from (4.14) we have 
$$D_{11\alpha\alpha}u = D_{11}uD_{\alpha\alpha}H  +  2 D_{11}u\left(D_{1\alpha}u \right)^2 +  { D_{11}u \over a}\left( D_{\alpha\alpha}u + \sum_{k= 1}^n a_kD_{\alpha\alpha k}u\right) , $$
which gives by replacing in (4.19) 
$$  D_{11}\partial_tu  \le - \left({1 \over F} - f\right) (D_{11}u)^2  + {D_{11}u \over F^2} \sum_{\alpha=1}^n  F_{\alpha\alpha}D_{\alpha\alpha}H -  {(D_{11}u)^2 \over F^2} \sum_{\alpha=1}^n  F_{\alpha\alpha}D_{\alpha\alpha}u$$
$$+  {D_{11}u \over a F^2} \sum_{\alpha =1}^n  F_{\alpha\alpha}D_{\alpha\alpha}u + {D_{11}u \over aF^2} \sum_{\alpha, k= 1}^n F_{\alpha\alpha} a_kD_{\alpha\alpha k}u  + D_{11}f  + D_{n+1}f D_{11}u , \eqno(4.20)$$
and since $F$ is homogenous of degree $k$ we have  at $(t_0, 0)$ 
$$\sum_{\alpha=1}^n  F_{\alpha\alpha}D_{\alpha\alpha}u = kF.$$
So  it follows from (4.20) that 
$$  D_{11}\partial_tu \le -  \left({k+ 1 \over F} - f\right)(D_{11}u)^2  +   k {D_{11}u \over a F} +  {D_{11}u \over F^2} \sum_{\alpha=1}^n  F_{\alpha\alpha}D_{\alpha\alpha}H $$
$$+ \  {D_{11}u \over aF^2} \sum_{\alpha, k= 1}^n  F_{\alpha\alpha}a_kD_{\alpha\alpha k}u  +   D_{11}f  + D_{n+1}f D_{11}u . \eqno (4.21) $$

\medskip

\noindent Since $H$ achieves a local maximum at $(t_0, 0)$, then the matrix $\left[D_{ij}H\right]$ is negative semi-definite at $(t_0, 0)$, and since $\left[F_{ij}\right]$ is positive semi-definite and diagonal at $(t_0 , 0)$, then we have at $(t_0, 0)$
$$ \sum_{\alpha=1}^n  F_{\alpha\alpha}D_{\alpha\alpha}H  \le 0.$$
Then  using the fact that $D_{11}u(t_0, 0) >0$,  we get  from (4.21) at $(t_0, 0)$,
$$  D_{11}\partial_tu  \le - \left({k+ 1 \over F} - f\right)(D_{11}u)^2   +   k {D_{11}u \over a F} +    {D_{11}u \over aF^2} \sum_{\alpha, k= 1}^n F_{\alpha\alpha}a_kD_{\alpha\alpha k}u  +  D_{11}f  + D_{n+1}f D_{11}u  . \eqno (4.22)$$
  
\medskip

\noindent Let us prove that  the first term in the right side of (4.22) is negative,that is 
$$ {k+ 1 \over F} - f \ge 0.  \eqno (4. 23) $$
 If $ 0< k \le 1$, then by Proposition 3.2, we have  a
$$ {1 \over F} - f   \ge 0 \eqno(4.24)$$
since $ {1 \over F} - f  = {\rho \over \sqrt{\rho^2 + |\nabla \rho|^2}  } \partial_t\rho \ge 0 $.  It is clear that (4.24) implies (4.23) since $F >0$.  Now if $ k > 1$, then by Proposition 3.2 we have 
$$ 0 \le - \left({1 \over F} - f\right)   = - {\rho \over \sqrt{\rho^2 + |\nabla \rho|^2}  } \partial_t\rho \le |\partial_t\rho|   \le  {R_2 \over R_1} \max_{x \in \mathbb{S}^n} \bigl|\mathcal{F}[\rho_0](x)\bigr| \  , $$
that is, 
$${1\over F} \ge f -  {R_2 \over R_1} \max_{x \in \mathbb{S}^n} \bigl|\mathcal{F}[\rho_0](x)\bigr| .  \eqno (4.25)$$
Now it is easy to see that (4.23) is a consequence of (4.25) and  the second part of  condition 4.1 in Proposition 4.1. Thus it follows from (4.22) and (4.23) that at $(t_0, 0)$,
$$ \partial_tD_{11}u  \le  k {D_{11}u \over a F} +    {D_{11}u \over aF^2} \sum_{\alpha, k= 1}^n F_{\alpha\alpha}a_kD_{\alpha\alpha k}u  +  D_{11}f  + D_{n+1}f D_{11}u  . \eqno (4.26)$$

\medskip

On the other hand, since at $(t_0, 0)$  we have 
$$D_ka_{ij} = D_{ijk}u,$$
then by differentiating equation (2.19) we get  at $(t_0, 0)$
$$D_k\partial_tu = {1 \over F^2}\sum_{i,j= 1}^nF_{ij}D_{ijk}u + D_kf   =   {1 \over F^2}\sum_{\alpha= 1}^nF_{\alpha\alpha}D_{\alpha\alpha k}u + D_kf  \eqno (4.27)$$
since $\left[F_{ij}\right]$ is diagonal at $(t_0, 0)$. 

\medskip

\noindent Now differentiating $H$ with respect to $t$, we see that at $(t_0, 0)$ 
$$\partial_tH = {\partial_tD_{11}u \over D_{11}u} - {\partial_t\varphi \over \varphi} =  {\partial_tD_{11}u \over D_{11}u}  + {1\over a} \partial_t u  - {1 \over a}\sum_{k=1}^na_kD_{k}\partial_tu $$
and using equation (2.18)  and (4.27) we obtain then 
$$\partial_tH= {\partial_tD_{11}u \over D_{11}u}  -  {1\over a} \left({1\over F} - f \right)  - {1 \over aF^2}\sum_{\alpha, k=1}^n F_{\alpha\alpha}a_kD_{\alpha\alpha k}u -{1\over a}\sum_{k=1}^na_kD_kf.   \eqno (4.28)$$
Thus we obtain from (4.26) and (4.28) 
$$ \partial_tH \le     {k - 1 \over a F} +   { D_{11}f  \over D_{11}u}  + D_{n+1}f   -{1\over a}\sum_{k=1}^na_kD_kf  +  {1\over a} f  \eqno (4.29)$$
at $(t_0, 0)$.   Since by (4.12) we have $\partial_tH(t_0, 0) \ge 0$, then it follows from (4.29) that 
$$ 0 \le    {k - 1 \over a F} +   { D_{11}f  \over D_{11}u}  + D_{n+1}f   -{1\over a}\sum_{k=1}^na_kD_kf  +  {1\over a} f  . \eqno (4.30)$$
And since
$$D_{n+1}f(a_1, .., a_n, -a)   -{1\over a}\sum_{k=1}^na_kD_kf(a_1, .., a_n, -a) = - {1\over a}\rho\partial_{\rho}f(a_1, .., a_n, -a), $$
then (4.30) becomes
$$ 0 \le    {k - 1 \over a F} +   { D_{11}f  \over D_{11}u}  - {1\over a}\rho\partial_{\rho}f +  {1\over a} f   \eqno (4.31)$$

\noindent  at $(t_0, 0)$. But by Proposition 3.2 we have $\partial_t\rho \ge 0$ if $k\le 1$, and $\partial_t\rho \le 0$ if $k > 1$. This implies that  
$\displaystyle {1 \over F(a_{ij})} - f(\rho x) \ge 0 $ if $k \le 1$, and $\displaystyle {1 \over F(a_{ij})} - f(\rho x) \le 0 $ if $k > 1$. That is, 
$${k-1 \over F(a_{ij})} \le (k-1) f(\rho x) . \eqno (4.32)$$

 \noindent It follows from  (4.31) and (4.32)    that 

$$ {1\over a}\left(  \rho \partial_{\rho}f  - k f\right)  \le { D_{11}f  \over D_{11}u} \eqno (4.33) $$

\medskip

\noindent at $(t_0 , 0)$.   Since $f$ satisfies  (1.7), then $ \rho\partial_{\rho}f- kf>0$, which implies
 $$ \delta_0 = \min_{(\rho , x) \in [ R_1  ,  R_2]\times\mathbb{S}^n} \bigl( \rho\partial_{\rho}f(\rho x)- kf(\rho x )  \bigr) >  0,$$
and  since $ R_1 \le \rho(t,x) \le R_2$ by Proposition 3.1, then   $ \rho \partial_{\rho} f   - kf  \ge \delta_0$. Thus we get from (4.33) at  $(t_0, 0)$

 $$ {\delta_0 \over a}  \le  { D_{11}f  \over D_{11}u} \le  { C  \over D_{11}u} \ ,  \eqno (4.34)$$
  where 
  $$C = \|f\|_{C^2( A_{R_1, R_2})} , \  \   \hbox{with} \   \  A_{R_1 , R_2} = \left\{ X \in \mathbb{R}^{n+1} \ : \  R_1 \le |X| \le R_2 \right\} . $$
 
 \noindent We recall that by definition of $H$, we have $D_{11}u = ae^{H}$ at $(t_0, 0)$. It follows from (4.34) that  
 
  $$ e^{H(t_0 , 0)}  \le  \delta_0^{-1}C  $$
or equivalently
 $$H(t_0, 0)  \le  \log{C \over \delta_0}  \ . \eqno (4.35) $$
  
  \noindent Thus the estimate $(4.3)$ is proved by taking 
  $$ C_0 =  \max\left( \log {C\over \delta_0} \  , \  \max_{x \in S^{n}}h(0, x)\right) . $$
  (4.3)  implies  then,  for any $(t, x) \in [0, T]\times \mathbb{S}^n$, 
  $$h(t, x) \le C_0 . \eqno (4.36) $$
  
 \noindent  We have by (4.2) 
  $$\max_{1 \le i \le n} \kappa_i = \langle X, \nu \rangle e^{h} $$
  and  since by Proposition 3.1 we have 
 $$  \langle X, \nu \rangle  = {\rho^2 \over \sqrt{\rho^2 + |\nabla\rho|^2}} \le \rho \le R_2 , $$
 then we get from (4.36) the upper bound 
 $$ \max_{1 \le i \le n} \kappa_i \le  R_2e^{C_0} .  \eqno (4.37)$$
 
 \noindent Now, to get a lower bound on the principal curvatures, it suffices to observe that $ \kappa_1 + \cdots + \kappa_n > 0$ since  $\kappa = (\kappa_1, ..., \kappa_n) \in \Gamma$, and  
 then use the upper bound (4.37).  Indeed, we have  for all $i = 1, ..., n$, 
 $$ 0 <  \kappa_1 + \cdots + \kappa_n \le \kappa_i + (n-1)R_2e^{C_0}$$
 so 
 $$\kappa_i \ge -(n-1)R_2 e^{C_0} .$$
 The proof of Proposition 4.1 is complete.

\end{proof}

\bigskip

The previous proposition allows  us  to get  higher order estimates on our solutions.

\begin{prop} Let $\rho : [0, T]\times\mathbb{S}^n  \to (0, +\infty)$ be an admissible  solution of (2.10) as in Proposition 4.1. Then for any $m \in \mathbb{N}$, there exist  two  positive constants $C_m$ and $\lambda_m$  depending only on $m, f, F , r_1, r_2$ and $\rho_0$ such that 
$$\|\rho  \|_{C^m([0, T] \times \mathbb{S}^n)}  \le C_m  \eqno (4.38) $$
and  for all $t \in [0, T]$, 
$$\|\partial_t\rho(t, .)\|_{C^m(\mathbb{S}^n)} \le C_me^{-\lambda_m t}  . \eqno (4.39) $$

\medskip

\noindent Moreover, there exists a compact set $K \subset  M(\Gamma)$  depending only on $f, F , r_1, r_2$ and $\rho_0$, sucht that  for any $ (t, x) \in [0, T] \times\mathbb{S}^n$,   
$$[a_{ij}(t, x)] \in K  \ ,     \eqno (4.40)$$
where the cone $M(\Gamma)$ is defined in section 2, and the matrix $[a_{ij}]$ is given by (2.4) in section 2. 

\end{prop}

\medskip

\begin{proof} The  principal curvatures $\kappa_i$ of the hypersurface $M_t$ parametrized by $X(t , x) = \rho(t, x)x$, are the eigenvalues of the matrix $[ a_{ij} ]$ (see section 2) defined by 
$$ [a_{ij}] = [ g^{ij} ]^{{1 \over 2}} [h_{ij}] [g^{ij}]^{{1 \over 2}} \eqno (4.41)$$
 where \ $ [g^{ij} ]^{{1 \over 2}}$ is the positive square root of \ $ [ g^{ij} ]$ which is given by
$$
 [ g^{ij} ]^{{1 \over 2}}= \rho^{-1} \left[ \delta_{ij} - {\nabla_i\rho\nabla_j \rho \over \sqrt{\rho^2 + | \nabla \rho |^2 }(\rho  + \sqrt {\rho^2 + | \nabla \rho |^2 })} \right]
 \eqno (4.42) $$
and $[h_{ij}]$ is the matrix representing the second fondamental form of $M_t$, given by 
$$
h_{ij} = \left ( \rho^2 + | \nabla \rho |^2  \right)^{-{1 \over 2}} (\rho^2\delta_{ij} + 2 \nabla_i\rho\nabla_j \rho - \rho\nabla_{ij}\rho) .  \eqno (4.43)$$

\noindent It is clear from Proposition 4.1, Proposition 3.1 and Proposition 3.3   by using (4.41), (4.42) and (4.43) that 
$$\sup_{t \in [0, T]} \|\rho(t, . ) \|_{C^2(\mathbb{S}^n)} \le C , \eqno (4.44) $$
where $C$ depends only on $f , r_1, r_2$ and $\rho_0$.  In order to get higher order estimates,  let us first prove  (4.40). By Proposition 3.2 we have 
$$|\partial_t \rho| \le Ce^{-\lambda t} \le C, \eqno (4.45)$$
where the constant  $C$ depends only on $f, r_1, r_2$ and $\rho_0$. Since $\rho$  satisfies (2.10), then it follows from (4.45) 
$${1 \over F(a_{ij})} - f(\rho x)  \le \left| {1 \over F(a_{ij})} - f(\rho x)\right|{ \sqrt{\rho^2 +  |\nabla\rho|^2 } \over \rho } = |\partial_t \rho| \le  C$$
that is,
$${1 \over F(a_{ij})}  \le f(\rho x) + C \le C_0$$
or equivalently 
$$F(a_{ij}) \ge {1 \over C_0} \  ,  \eqno (4.46)$$
where 
$$C_0 = C + \max_{R_1 \le |X| \le R_2}|f(X)| \ .$$
 Since $F \equiv 0$ on $\partial M(\Gamma)$,   it follows from (4.46) that there exists a constant $  \delta_0 > 0$ depending only on  $f, F, R_1, R_2$ and $\rho_0$ such that 
$$\text{dist}\Bigl([a_{ij}], \partial M(\Gamma)\Bigr) \ge \delta_0 , \eqno (4.47) $$
where $\partial M(\Gamma)$ is the boundary of the cone $M(\Gamma)$ and $\text{dist}\bigl([a_{ij}] ,  \partial M(\Gamma)\bigr) $ is the distance of $[a_{ij}]$ to $\partial M(\Gamma)$. It is clear from (4.47) that there exists a compact set $K \subset M(\Gamma)$ depending only on $f, F, r_1, r_2$ and $\rho_0$ such that  \ $ [a_{ij}] \in K $. Thus $(4.40)$ is proved.

Let us now prove the estimates (4.38) and (4.39). Since $F$ satisfies (1.3)(or equivalently (2.6)), it follows from (4.40)  and the estimate (4.44) that  equation (2.10) is uniformly parabolic. Since by hypothesis the function $\log F$ is concave, then we can apply  a result of B. Andrews \cite{bA04} (Theorem 6, p.3 ), which is a generalisation of the result of N. Krylov  \cite{nK87} on fully nonlinear parabolic equations,  to obtain the estimate  
$$ \|\partial_t\rho \|_{C^{\alpha}([0, T]\times \mathbb{S}^n)} +   \|\nabla_{ij}\rho \|_{C^{\alpha}([0, T]\times \mathbb{S}^n)} \le C , \eqno (4.48)$$
where $C^{\alpha}([0, T]\times \mathbb{S}^n)$ is the parabolic H\"older's space, and 
where the constants $C>0$, $\alpha \in (0, 1)$ depend only on $f, F, r_1, r_2 $ and $\rho_0$.  The higher order estimates (4.38) follows from (4.48) and the standard theory of linear   parabolic equations (see \cite{oL}).   In order to prove (4.39) we use the following well known interpolation inequality, which is valid on any compact Riemannian manifold $M$, 
$$\|\nabla u\|_{L^{\infty}(M)}^2 \le 4 \|u\|_{L^{\infty}(M) }  \|\nabla^2 u\|_{L^{\infty}(M)}, \  \    u \in C^{\infty}(M) ,  \eqno (4.49)$$
where $\nabla u$ and $\nabla^2u$ denote respectively  the gradient and the hessian  of $u$. It suffices to apply (4.49) first to $u= \partial_t\rho$ and iterate it on the spatial higher order derivatives of $\partial_t \rho$  and using (4.38) and (3.11)  to get (4.39). This achieves the proof of Proposition 4.2.

\end{proof}
\bigskip

Now we are in position to prove our main result.

\medskip

\begin{proof}[Proof of Theorem 1.1 and Theorem 1.2 ]  Let $X_0(x) = \rho_0(x)x$ satisfies conditions (1.10) in Theorem 1.1 or conditions (1.13) in Theorem 1.2.  Let  $X : [0, T] \times \mathbb{S}^n \to \mathbb{R}^{n+1}$ a local solution of (1.2).   As we saw in section 2,  $X$ is given by 
$$X(t, x) = \rho(t, \varphi(t, x))\varphi(t, x) , \  (t, x) \in [0 , T]\times\mathbb{S}^n \eqno (4.50) $$
where $\rho$ satisfies (2.10) and $\varphi( t , .) : \mathbb{S}^n \to \mathbb{S}^n$ is a diffeomorphism satisfying the ODE 
$$\begin{cases} \partial_t\varphi(t, x) = Z(t, \varphi(t, x)  )
 \cr \varphi(0 , x) = x  ,  \end{cases} \eqno (4.51)$$
 with 
 $$Z(t, y) =  -  \left( {1 \over F(a_{ij}(t,y))}  - f(\rho(t, y) y) \right) {\nabla \rho(t,y) \over \rho  \sqrt{|\nabla \rho(t,y)|^2 + \rho^2(t, y)}} $$
 
       $$           = -  {\partial_t\rho(t , y)\nabla\rho(t,y) \over |\nabla \rho(t,y)|^2 + \rho^2(t, y) } \ 
  ,  \  (t, y) \in [0, T]\times\mathbb{S}^n .      \eqno (4.52) $$

\medskip

Since $X_0$ satisfies condition (1.10) in Theorem 1.1 or condition (1.13) in Theorem 1.2,   then it is easy to check that the hypothesis of Proposition 4.1(and then Proposition 4.2)   concerning $\rho_0$ are  satisfied.  We can then apply Proposition 4.2  to the function  $\rho$ given above.  If we differentiate equation (4.51) and use the estimates (4.38)-(4.39) in Proposition 4.2, then it is not difficult to see that for any $m \in \mathbb{N}$, we have 
$$\|\varphi\|_{C^m([0, T]\times\mathbb{S}^n , \  \mathbb{S}^n)} \le C_m \eqno (4.53) $$
and  for any $t \in [0, T]$, 
$$\|\partial_t\varphi(t , . )\|_{C^m(\mathbb{S}^n,\ \mathbb{R}^{n+1})} \le C_me^{-\lambda_m t } ,  \eqno (4.54) $$
where $C_m$ and $\lambda_m$  are  positive constants depending only on $m, f, F, r_1, r_2$ and $X_0$. It follows from (4.50) by using  the estimates (4.38)-(4.39)  in Proposition 4.2 and (4.53)-(4.54)  that, for any $m \in \mathbb{N}$, 
$$\|X\|_{C^m([0, T]\times\mathbb{S}^n , \ \mathbb{R}^{n+1})} \le C_m \eqno (4.55) $$
and  for all $t \in [0, T]$, 
$$\|\partial_tX(t, .)\|_{C^m(\mathbb{S}^n , \ \mathbb{R}^{n+1})} \le C_me^{-\lambda_m t}, \eqno (4.56) $$ 
with new constants  $C_m$ and $\lambda_m$   depending only on $m, f, F, r_1, r_2$ and $X_0$. Also by Proposition 4.2 there exists a compact set $K \subset M(\Gamma)$ depending only on  $m, f, F, r_1, r_2$ and $X_0$ such that  for any $(t, x) \in [0, T]\times\mathbb{S}^n$, we have 
$$[a_{ij}(t, x)] \in K \subset M(\Gamma)  \ ,     \eqno (4.57)$$
where the matrix $[a_{ij}]$ is given by (2.4).  Since the constant $C_m$   in (4.55) and the compact set $K$  in (4.57)  are independant of $T$, then $X$  can be extended to $[0, +\infty)$  as a solution of (1.2). The estimates (4.55), (4.56)  and  (4.57)  become then 
$$ \|X \|_{C^m([0, + \infty) \times \mathbb{S}^n,\  \mathbb{R}^{n+1})}  \le C_m  \eqno (4.58)$$
 
$$\|\partial_tX(t, .)\|_{C^m(\mathbb{S}^n, \ \mathbb{R}^{n+1})} \le C_me^{-\lambda_m t} \  \   \hbox{for all} \   t \in [0, + \infty) \eqno (4.59) $$ 
and
$$[a_{ij}(t, x)] \in K  \subset M(\Gamma)  \  \   \hbox{for all} \   t \in [0, + \infty) .   \eqno (4.60)$$

\medskip

Now it is clear from (4.58) and (4.59) that there exists a map $ X_{\infty} \in C^{\infty}(\mathbb{S}^n , \mathbb{R}^{n+1}) $ such that $X(t , .) \to X_{\infty}$ as $t \to + \infty$ in 
$C^m(\mathbb{S}^n , \mathbb{R}^{n+1})$ for all $ m\in \mathbb{N}$, and satisfying 
$$\| X(t , .) - X_{\infty}\|_{C^m(\mathbb{S}^n, \ \mathbb{R}^{n+1})} \le C_me^{-\lambda_m t} \  \   \hbox{for all} \   t \in [0, + \infty) . $$

\medskip

 \noindent Since $X(t, .)$ is starshaped, then it is easy to see that $X_{\infty}$ is also starshaped, and from (4.60) we deduce that the principal curvatures of $X_{\infty}$  lie in $\Gamma$. 
By passing to the limit in equation (1.2)  and using (4.59), we see that $X_{\infty}$ satisfies 
$${1 \over F(\kappa(X_{\infty}) )} - f(X_{\infty}) = 0 .$$

 \noindent This achieves the proofs of Theorem 1.1 and Theorem 1.2.

\end{proof}

\begin{proof}[Proof of Corollary 1.1 ]  As in Remark 1.1, if we take $X_0(x) = rx $, where  $0 < r \le r_1$ with $r_1$ as in (1.8), then by using (1.7) and (1.8) one easily  checks that condition (1.10) in Theorem 1.1 is satisfied by $X_0$. Thus the evolution problem (1.2) admits a  global solution $X(t, .)$ which  converges as $t \to + \infty$,  to a solution $X_{\infty}$ of  
$$ {1 \over F(\kappa(X_{\infty}))} = f(X_{\infty})  \eqno (4.61)$$ 
which is starshaped and satisfying $\kappa(X_{\infty}) \in \Gamma$. It remains then to prove that $X_{\infty}$ is the unique starshaped solution of (4.61) such that $\kappa(X_{\infty}) \in \Gamma$.  Let $X_1$ and $X_2$ two   starshaped solutions of (4.61) such that  
$\kappa(X_l) \in \Gamma , \ l=1, 2$. We have then 
$${1 \over F(\kappa(X_l))} = f(X_l) \ , \  l = 1, 2.  \eqno (4.62)$$

Let $\rho_l  \ ( l= 1, 2 ) $  be the radial function of $X_l$, and set $u_l(x) = \log\rho_l(x)$.  Then we have by using formula (2.4) of section 2, 
$${ 1 \over F(a_{ij}(u_l)) }  = f(e^{u_l} x) \ , \  l = 1, 2,  \eqno (4.63) $$
where the matrix $[a_{ij}(u_l)] $ is given by
$$[ a_{ij}(u_l)] = {e^{-u_l} \over \sqrt{1 + |\nabla u_l|^2}}[\gamma_{ij}][b_{ij}][\gamma_{ij}]   \eqno (4.64) $$
with
$$\begin{cases}  \displaystyle b_{ij} = \delta_{ij} + \nabla_{i}u_l\nabla_j u_l -\nabla_{ij}u_l   \cr 
\cr \displaystyle 
\gamma_{ij} = \delta_{ij} -{\nabla_iu_l\nabla_ju_l \over\sqrt{1 + |\nabla u_l|^2} \left(1 + \sqrt{1 + |\nabla u_l|^2}\right) } \  , \  l= 1, 2.  \end{cases} \eqno (4.65) $$

\medskip

\medskip

We shall prove that  for any $x \in \mathbb{S}^n$, we have 
$$u_1(x)  \ge u_2(x) .  \eqno (4.66) $$
It is clear that (4.66) would imply that $u_1 = u_2$, and then $\rho_1 = \rho_2$.  To prove  (4.66)  define a function $u :  \mathbb{S}^n \to \mathbb{R}$ by $u(x) = u_1(x) - u_2(x)$, and   let $x_0 \in \mathbb{S}^n$ a  point where $u$ achieves its minimum. Then  we have at  $x_0$ that  $ \nabla u = 0 $ and the matrix $\nabla^2u$ is positive semi-definite, that is,   
$ \nabla u_1 = \nabla u_2 $ and  $\nabla^2u_1 \ge \nabla^2u_2 $ (in the sense of operators) at $x_0$.  This implies by using     (4.64) and (4.65)  that  at $x_0$, 
$$e^{u_1}[a_{ij}(u_1)] \le e^{u_2} [a_{ij}(u_2)]  \eqno (4.67)$$
in the sense of operators.  Since the function $F$ is monotone (by (1.3))  and homogenous of degree $k$,   it follows from (4.63)  and (4.67) that 
$$e^{- ku_1(x_0)}f(e^{u_1(x_0)}x_0)  \ge  e^{- ku_2(x_0)}f(e^{u_2(x_0)}x_0)  $$
which implies by using (1.7) that $u_1(x_0)  \ge u_2(x_0)$ or equivalently $u(x_0) \ge 0$. This proves (4.66) and the proof of Corollary 1.1 is complete.

\end{proof}

\end{document}